\documentclass[12pt,twoside]{amsart}

\newtheorem{Thm}{Theorem}
\newtheorem{Cor}[Thm]{Corollary}
\newtheorem{Lemma}[Thm]{Lemma}
\newtheorem{Prop}[Thm]{Proposition}
\newtheorem*{Prop*}{Proposition}

\theoremstyle{definition}
\newtheorem{Defn}[Thm]{Definition}
\newtheorem{Ex}{Example}

\input{diagrams}

\title{The noncommutative geometry of the discrete Heisenberg group}
\author[T.~Hadfield]{Tom Hadfield}
\address{Department of Mathematics, University of California, Berkeley, CA 94720}
\email{hadfield@math.berkeley.edu}
\date{\today}

\begin{document}

\begin{abstract}
 Motivated by the search for new examples of ``noncommutative manifolds'', we study the noncommutative geometry of the group C*-algebra of the three dimensional discrete Heisenberg group. 
 We present a unified treatment of the K-homology, cyclic cohomology and derivations of this algebra. 
\end{abstract}

\maketitle

\section{The discrete Heisenberg group}

Recently, there has been considerable interest in the notion of a ``noncommutative manifold''. 
There are several different approaches to this problem. 
 The first is the operator-algebraic approach pioneered by Connes in [Co96], 
 which takes as its starting point the commutative algebra ${C^{\infty}}(M)$, where $M$ is a compact Riemannian spin manifold, together with the Dirac operator coming from the spin connection, and then gives an axiomatic formulation which extends to noncommutative algebras. 
 The most accessible example of a ``noncommutative differentiable manifold'' in this sense is provided by the noncommutative tori [Ri81].
 Recently Connes' program has enjoyed considerable success [CD01], [CL00], [Var01].  

A second approach, via quantum groups, has been much more examples driven, proceeding with  case by case studies of the many interesting algebras arising from quantum groups
[Ma00], [Sch99], [Wo87].
 Rather than constructing spectral triples, these workers focus on classifying the possible differential calculi over a given algebra. 
 
The route we take in this paper falls between these two approaches. 
 We start with a specific example, and see how much of Connes' formalism is applicable. 
   In Connes' picture, a noncommutative geometry over a C*-algebra $A$ consists of a spectral triple over $A$ equipped with additional structures. 
 The bounded formulation of spectral triples are 
 Fredholm modules, equivalence classes of which make up
the Kasparov K-homology groups ${KK^i}(A,{\bf C})$ ($i=0,1$). 
 Thus the focus of this work is to calculate the K-homology groups and exhibit the generating Fredholm modules. 
 We will not at this point address the problem of finding corresponding spectral triples, and then trying to equip these with the additional structures of a noncommutative geometry.

 We will study in detail the group C*-algebra of the discrete Heisenberg group $H_3$. 
This algebra is closely related to the noncommutative tori $A_{\theta}$  - 
it arises as a continuous field of noncommutative tori over the circle - and many of the constructions used for the $A_{\theta}$ can also be applied here. 
However, not everything generalizes so nicely. 
The group C*-algebra $C^{*}(H_3)$ is a crossed product $C({\bf T}^2) \times_{\alpha} {\bf Z}$, and provided the main motivation for our work on K-homology of crossed products by ${\bf Z}$ [Ha01].

The discrete three-dimensional  Heisenberg group $H_3$ 
can be defined abstractly as the group generated by elements 
$a$ and 
$b$ 
such that the commutator 
$c = ab {a^{-1}}{b^{-1}}$ 
is central. 
It  
can be realised as the multiplicative group of upper-triangular matrices
\begin{equation}
H_3 = 
\{
\left(
\begin{array}{ccc}
1 & a & c \cr
0 & 1 & b \cr
0 & 0 & 1 
\end{array}
\right)
 : a,b,c \in {\bf Z} \}.
\end{equation}
 It is a semidirect product 
$ H_3 \cong  { {\bf Z}^2} {\times}_{\alpha}  {\bf Z}$, 
with the action 
$\alpha$ of ${\bf Z}$ on 
${\bf Z}^2$ being given by 
${\alpha^k}(m,n)=(m, km+n)$.
Hence 
${H_3}$ 
is amenable.  
It is the prototypical example of a two-step nilpotent non-type I group.

Let 
$A= {C^{*}}( {H_3} )$. 
Since 
$H_3$ 
is amenable, 
${C^{*}}( {H_3}) \cong {C^{*}_{r}} ( {H_3})$. 
The group C*-algebra $A$ 
is the universal C*-algebra generated by unitaries 
$U$, $V$ and  $W$ 
satisfying 
$VU= WUV$, 
with $W$ central. 
 It is isomorphic to the crossed product algebra 
$C( {\bf T}^2 ) \times_{\alpha} {\bf Z}$,
where the action 
$\alpha$ 
of 
${\bf Z}$ 
on 
 ${C^{*}}(U,W) \cong$
$C({\bf T}^2)$ 
is given by
 $VU{V^*} = U^*$, $VW{V^*}=W$. 
 Given an irreducible representation 
$\pi : {C^{*}}(H_3) \rightarrow {\bf B}(H)$, 
 since $W$ is central, we must have $\pi(W) = \lambda I$, for some $\lambda \in {\bf C}$, $| \lambda | =1$, $\lambda = \exp (2 \pi i \theta)$. 
For each given 
$\theta$, 
we have a surjective *-homomorphism
$ \phi : C^{*}(H_3) \rightarrow {A_{\theta}} $
given by 
$ U \mapsto U$, 
$V \mapsto V$, 
$W \mapsto \lambda I$.
Hence ${C^{*}}(H_3)$ is a continuous field of rotation algebras $A_{\theta}$ over the circle. 

 We will calculate the K-homology of the group C*-algebra, and exhibit generating Fredholm modules. 
We compute the pairings of these Fredholm modules with the generators of K-theory, 
and investigate their behaviour in the six term cyclic sequence on K-homology dual to the corresponding Pimsner-Voiculescu sequence for K-theory.  
Having done this, 
we study the derivations of the group ring ${\bf C}H_3$, and of a dense smooth subalgebra 
$A^{\infty}$ of the group C*-algebra. 
Then we turn our attention to cyclic cohomology.
We apply Burgheleas' theorem to calculate the cyclic cohomology of the group ring 
${\bf C}{H_3}$.

\section{K-theory}

\begin{Thm}[Anderson-Paschke]
\label{kthyhbg}
[AP89]
The K-groups of $A$ are both 
${\bf Z}^3 $.
 The group 
${K_0}(A)$ 
is generated by the equivalence classes of  projections 
$1 \in A$ 
and 
${P_{a}}$, ${P_{b}} \in {M_2}(A)$.
Generators of 
${K_1}(A)$ 
correspond to the  unitaries 
$U, V \in A$ 
and 
${V_{a}} \in {M_2}(A)$.
\end{Thm}

The matrices  $P_a$, $P_b$ and  $V_a$ are defined in [AP89].
 The projections $P_a$ and $P_b$ correpond to the Bott projections for the commutative *-subalgebras 
${C^{*}}(U,W)$ and ${C^{*}}(V,W)$, 
 which are both isomorphic to $C({\bf T}^2)$.
 The unitary matrix $V_a$ is given by 
\begin{equation}
\label{Va}
V_a = Z_a ( V \otimes I_2 ) P_a + (I_2 - P_a)
\end{equation}
 where $Z_a$ is a unitary matrix satisfying $\alpha(P_a) =$
 $Z_a P_a {Z_a^{*}}$.

 For future reference, we summarize the behaviour of these generators  
 under  the Pimsner-\\ 
 Voiculescu six term cyclic sequence [PV80] 
 for K-theory of a crossed product by ${\bf Z}$. 
These results are taken from [AP89]. 
Taking $C({\bf T}^2)$ to be generated by the unitaries $U$ and $W$, with 
the automorphism $\alpha$ acting via $\alpha(U) = WU$, $\alpha(W) = W$, and  $C({\bf T}^2) \times_{\alpha} {\bf Z} \cong {C^{*}}(H_3)$, we have :  

\begin{diagram}
{K_0} ( C( {\bf T}^2 )) & \rTo^{ id - {\alpha_{*}}} & {K_0} (C( {\bf T}^2 )) & \rTo^{ i_{*}} 
& {K_0} ( {C^{*}}( H_3)) \\
\uTo^{\delta_1} &    &  &  & \dTo^{\delta_0} \\
{K_1} ({C^{*}}(H_3)) & \lTo^{ i_{*}} & {K_1} (C( {\bf T}^2 )) & \lTo^{id - {\alpha_{*}}} & {K_1} (C( {\bf T}^2 )) 
\end{diagram}\\
 Starting in the top left corner, 
${K_0}(C({\bf T}^2))$ is generated by $[1]$ and $[P_a]$. 
The map 
\begin{equation}
(id - {\alpha_{*}}) : {K_0}(C({\bf T}^2)) \rightarrow {K_0}(C({\bf T}^2))
\end{equation} 
is the zero map. 
For the boundary map $\delta_0$, we have 
\begin{equation}
\label{deltazero}
\delta_0 [1] = 0, \quad 
\delta_0 [P_a] = 0, \quad
\delta_0 [P_b] = [W].
\end{equation}
 The group ${K_1}(C({\bf T}^2))$ is generated by $[U]$ and $[W]$. 
 Under the map 
\begin{equation}
(id - {\alpha_{*}}) : {K_1}(C({\bf T}^2)) \rightarrow {K_1}(C({\bf T}^2))
\end{equation} 
 we have 
$(id - {\alpha_{*}})(W) = 0$ and  $(id - {\alpha_{*}})(U) = W$. 
 So ${i_{*}}[W] = 0$, and ${i_{*}}[U] = [U]$.
 For the boundary map $\delta_1$, we have 
\begin{equation}
\label{deltaone}
\delta_1 [U] = 0, \quad
\delta_1 [V] = [1], \quad
\delta_1 [V_a] = [P_a].
\end{equation}

\section{Fredholm modules as K-homology}

We begin with some preliminaries about Fredholm modules. 
Recall from [Co94], p288 that a Fredholm module over a *-algebra 
$A$ 
is a triple 
$(H,\pi,F)$,
where $\pi$ is a *-representation of 
$A$ 
as bounded operators on the Hilbert space 
$H$. 
The operator 
$F$ is a selfadjoint element of  
${\bf B}(H)$,
satisfying $F^2 =1$, 
such that the commutators 
$[F,\pi(a)]$ are compact operators for all 
$a \in A$. 
Such a Fredholm module is called odd. 
An even Fredholm module is the above data, 
together with a 
${\bf Z}_2$-grading 
of the Hilbert space $H$,
given by a grading operator 
$\gamma \in {\bf B}(H)$
with 
$\gamma = {\gamma}^{*}$, 
$\gamma^2 =1$,
$[\gamma, \pi(a)]=0$ 
for all 
$a \in A$, 
and
$F\gamma = - \gamma F$. 
 The *-algebra 
$A$ 
will usually be a dense subalgebra of a C*-algebra, 
closed under holomorphic functional calculus.

We define simple even and odd Fredholm modules that we will use extensively throughout this work.

\begin{Ex}
\label{evenfred}
Given a C*-algebra $A$, with a *-homomorphism $\phi : A \rightarrow {\bf C}$, 
we construct a canonical even Fredholm module
 ${\bf z}_0 \in {KK^0}(A, {\bf C})$ : 
\begin{equation}
{\bf z}_0 = (
H_0 = {\bf C}^2, \pi_0 = \phi \oplus 0 , 
F_0 = \left(
\begin{array}{cc}
0 & 1 \cr
1 & 0 \cr
\end{array}
\right),
\gamma = 
\left(
\begin{array}{cc}
1 & 0 \cr
0 & -1 \cr
\end{array}
\right)
).
\end{equation}
 In general ${\bf z}_0$ may well represent a trivial element of the even K-homology of $A$ (for example, if $\phi$ is the zero homomorphism.) 
However,

\begin{Lemma} 
\label{evenfrednontriv}
If $A$ is unital, and $\phi$ nonzero, then 
 $< {ch}_{*}( {\bf z}_0 ), [1]> =1$.
Hence in this situation ${\bf z}_0$ is a  nontrivial element of K-homology.
\end{Lemma}
\begin{proof}
Here, 
${ch}_{*} : {KK^0}(A,{\bf C}) \rightarrow {HC^{even}}(A)$, 
is the even Chern character as defined in [Co94], p295, mapping the even K-homology of $A$ into even periodic cyclic cohomology,
and $<.,.>$ denotes the pairing between K-theory and periodic cyclic cohomology defined in [Co94], p224.
 We have
\begin{equation}
 < {ch_{*}}( {\bf z}_0 ), [1] > = 
{{\lim}_{n \rightarrow \infty} } {(n!)}^{-1} \psi_{2n}  (1,...,1),
\end{equation}
 where (for each $n$)  $\psi_{2n}$ is the cyclic $2n$-cocycle defined by
\begin{equation}
 \psi_{2n} ( a_0, a_1, ..., a_{2n}) = 
(-1)^{n(2n-1)} \Gamma (n+1) 
 Tr( \gamma \pi_0 (a_0) [ F_0 , \pi_0 (a_1)] ... [F_0 , \pi_0 (a_{2n})]). 
\end{equation}
 Since $\Gamma(n+1) = n!$ it follows that
\begin{equation}
< {ch_{*}}( {\bf z}_0 ), [1] > = 
{{\lim}_{n \rightarrow \infty} }
 {(-1)}^n
 Tr( \gamma \pi_0 (1) {[F_0 , \pi_0 (1)]}^{2n}).
\end{equation} 
 Now,  
$[F_0 , \pi_0 (1) ] = 
\left(
\begin{array}{cc}
0 & -1 \cr
1 & 0 \cr
\end{array}
\right),$
 hence
 $\gamma \pi_0 (1) {[F_0 , \pi_0 (1)]}^{2n} =$
$(-1)^n 
\left(
\begin{array}{cc}
1 & 0 \cr
0 & 0 \cr
\end{array}
\right)$.
 Therefore
\begin{equation} 
< {ch_{*}}( {\bf z}_0 ), [1] > = 
 {{\lim}_{n \rightarrow \infty} } {(-1)}^n 
Tr( 
(-1)^n 
\left(
\begin{array}{cc}
1 & 0 \cr
0 & 0 \cr
\end{array}
\right)) = 1
\end{equation}
as claimed.
 It also follows that $[1]$ is a nontrivial element of $K_0 (A)$. 
\end{proof}
\end{Ex}

\begin{Ex}
\label{oddfred}
 Let $A$ be a C*-algebra, together with a *-homomorphism 
$\phi : A \rightarrow {\bf C}$ and an automorphism $\alpha$ implementing an action of ${\bf Z}$. 
We describe a canonical odd Fredholm module
 ${\bf z}_1 \in {KK^1}( A \times_{\alpha} {\bf Z}, {\bf C})$ : 
\begin{equation}
 {\bf z}_1 = (H_1 ={l^2}( {\bf Z}), \pi_1, F_1)
\end{equation}
 Take 
$\pi_1 : A \times_{\alpha} {\bf Z} \rightarrow {\bf B}( {l^2}({\bf Z}))$ 
to be defined by 
\begin{equation}
(\pi_1 (a) \xi) (n) = \phi( {\alpha^{-n}}(a)) \xi(n), \quad
(\pi_1 (V) \xi) (n) = \xi(n-1),
 \end{equation} 
 for $\xi \in {l^2}({\bf Z})$, $a \in A$, and $V$ the unitary implementing the action of ${\bf Z}$ on $A$ 
(via $Va{V^{*}} = \alpha (a)$).  
 Then $\pi_1$ is the usual representation of $A \times_{\alpha} {\bf Z}$ induced from the representation $\phi$ of $A$.
 We take 
\begin{equation}
F_1  \xi(n) = sign(n) \xi(n) =
\left\{
\begin{array}{cc}
 \xi(n) & : n \geq 0 \cr
- \xi(n) & : n<0 \cr
\end{array}
\right. 
\end{equation}
 The Fredholm module  ${\bf z}_1$ is always  nontrivial
(even if $\phi$ is the zero homomorphism). We have:

\begin{Lemma}
$< {ch}_{*} ({\bf z}_1 ), [V]> =1$. 
\end{Lemma}
\begin{proof}
Again, 
${ch}_{*} : {KK^1}(A,{\bf C}) \rightarrow {HC^{odd}}(A)$, 
is the odd Chern character as defined in [Co94], p296, mapping the odd K-homology of $A$ into odd periodic cyclic cohomology,
and $<.,.>$ denotes the pairing between K-theory and periodic cyclic cohomology [Co94], p224.
 Instead of calculating the pairing directly, as in the previous example, we instead use Connes' index theorem [Co94], p296, which states that 
\begin{equation}
< {ch_{*}}( {\bf z}_1 ), [V]> 
= Index (EVE)
\end{equation}
where $E = {\frac{1}{2}}(1+F)$ is the natural orthogonal projection 
 ${l^2}({\bf Z}) \rightarrow {l^2}({\bf N})$. We have 
\begin{equation}
Index(EVE) = dim\, ker(EVE) - dim\, ker(E{V^{*}}E) = 1 -0 =1,
\end{equation}
hence the result.
 This shows that ${\bf z}_1$ is a nontrivial Fredholm module, 
 and also that $[V] \neq 0 \in$ 
${K_1}(A \times_{\alpha} {\bf Z})$.
 \end{proof}
\end{Ex}

\section{K-homology}

In this section and the next we study the K-homology of the group C*-algebra ${C^{*}}(H_3)$. 
We exhibit the generating Fredholm modules, and study their behaviour in the six term cyclic sequence on K-homology for crossed products by ${\bf Z}$. 
 Since the discussion is both lengthy and  technical, we note for the convenience of the reader that the most important results are summarized in Theorems 19 and 23.
 We begin with:

\begin{Prop}
\label{khomhbg}
\label{kgpsfreeab}
The K-homology of ${C^{*}}(H_3)$ is given by
 ${KK^i}( {C^{*}}(H_3) , {\bf C}) \cong {\bf Z}^3$ 
$(i=0,1)$. 
\end{Prop}
\begin{proof} 
 It is an immediate consequence of Rosenberg and Schochet's universal coefficient theorem [RS87], 
that for a separable C*-algebra $A$, belonging to the ``bootstrap class'', if the K-groups ${K_i}(A)$ ($i=0,1$) are free abelian, 
then we have  isomorphisms ${KK^i}(A,{\bf C}) \cong$ ${K_i}(A)$ (as abelian groups). See also [Bla98], p234.
 Hence for ${C^{*}}(H_3)$, since the K-groups are free abelian,
 the result is immediate.
\end{proof}

Now we want to identify the generating Fredholm modules. 
We start by considering the noncommutative tori.
  Recall that 
${C^{*}}(H_3)$ is a continuous field of rotation algebras $A_{\theta}$ over the circle. For each 
$\theta \in {\bf R}$, 
we have a surjective *-homomorphism
$ \phi : C^{*}(H_3) \rightarrow {A_{\theta}} $
given by 
$ U \mapsto U$, 
$V \mapsto V$, 
$W \mapsto \lambda I$, where $\lambda = \exp (2 \pi i \theta)$.
 We can pull back known  Fredholm modules over the 
$A_{\theta}$ via $\phi$  
to give Fredholm modules over 
$C^{*} (H_3)$ : 
\begin{equation}
{\phi}^{*} :  {KK^{i}}( {A_{\theta}}, {\bf C} ) \rightarrow  {KK^{i}} ( {C^{*}}( H_3 ) , {\bf C} ).
\end{equation}
It will be sufficient for our purposes just to do this for the case $\theta = 0$, and henceforth $\phi$ will denote the *-homomorphism 
$\phi : {C^{*}}( H_3 ) \rightarrow C( {\bf T}^2 )$, 
with $\phi(W) = I$. 
 We summarize known results about the K-homology of $C({\bf T}^2)$.

\begin{Prop} 
\label{khomCT2}
The even K-homology ${KK^0}( C({\bf T}^2), {\bf C}) \cong {\bf Z}^2$ is generated by Fredholm modules ${\bf w}_0$ and ${\bf Dirac}$.
 The odd K-homology ${KK^1}( C({\bf T}^2), {\bf C}) \cong {\bf Z}^2$ is generated by Fredhom modules ${\bf w}_1$ and ${\bf w}_1 '$.
We give an explicit description of each of these Fredholm modules. 
\end{Prop} 
\begin{proof} 
 Since the K-groups ${K_0}(C({\bf T}^2))$ and ${K_1}(C({\bf T}^2))$
 are both isomorphic to ${\bf Z}^2$ it again follows from the universal coefficient theorem (Prop 4) that both the even and the odd K-homology groups are isomorphic to ${\bf Z}^2$ also.
 We will consider $C({\bf T}^2)$ to be the abelian  C*-algebra generated by commuting unitaries $U$ and $V$. 
The odd K-homology of 
$C( {\bf T}^2 )$ is as follows.
 We can think of $C({\bf T}^2)$ as a crossed product 
$C({\bf T}) \times {\bf Z}$ via a trivial action of ${\bf Z}$. 
 Taking $C({\bf T}) \cong$ ${C^{*}}(V)$, and then 
$C({\bf T}) \cong$ ${C^{*}}(U)$, with the trivial ${\bf Z}$-action being implemented by $U$ and then $V$ respectively, Example 2 gives us generating Fredholm modules 
${\bf w}_1 = ( {l^2}({\bf Z}),{\pi_1},F)$ and 
${{\bf w}_1}'=( {l^2}({\bf Z}), {\pi_1}' , F)$, 
with $\pi_1$ and ${\pi_1}'$ given by 
\begin{equation}
\label{wone}
{\pi_1}(U) = S,\quad 
{\pi_1}(V) = I, \quad
{\pi_1}'(U) = I,\quad
{\pi_1}'(V) = S.
\end{equation}
Here $S$ denotes the shift $S{e_n} = e_{n+1}$ with respect to the canonical orthonormal basis 
${\{ {e_{n}} \}}_{n \in {\bf Z}}$ of 
 ${l^2}({\bf Z})$.
In each case $F$ is the diagonal operator
$F{e_n} = sign(n){e_n}$.

We describe the generators ${\bf w}_0$ and ${\bf Dirac}$ of the even K-homology. 
 The Fredholm module ${\bf w}_0$ is the canonical even Fredholm module (Example 1) corresponding to the unital *-homomorphism 
$\phi : C({\bf T}^2) \rightarrow {\bf C}$ 
given by 
$U, V \mapsto 1$ :
\begin{equation}
{\bf w}_0 = (
H = {\bf C}^2, \pi = \phi \oplus 0 , 
F = \left(
\begin{array}{cc}
0 & 1 \cr
1 & 0 \cr
\end{array}
\right),
\gamma = 
\left(
\begin{array}{cc}
1 & 0 \cr
0 & -1 \cr
\end{array}
\right)
).
\end{equation} 
 The even Fredholm module 
${\bf Dirac}$ is 
the bounded formulation of the Dirac operator on ${\bf T}^2$, 
\begin{equation}
\label{T2dirac}
{\bf Dirac} = (H, \pi, F)
\end{equation} 
 where $H = {l^2}( {\bf Z}^2) \oplus {l^2}( {\bf Z}^2)$. 
We have $C( {\bf T}^2)$ acting on ${l^2}({\bf Z}^2)$ via 
\begin{equation}
 \pi_0 (U){e_{m,n}} = e_{m+1,n}, \quad 
\pi_0 (V){e_{m,n}} = e_{m,n+1}
\end{equation} 
 and we take
\begin{equation}
\label{DiracF}
 \pi(a)=
\left(
\begin{array}{cc}
\pi_0 (a) & 0 \cr
0 & \pi_0 (a) \cr
\end{array}
\right), \quad
 F=
\left(
\begin{array}{cc}
0 & {F_0} \cr
{F_0}^{*} & 0 \cr
\end{array}
\right)
\end{equation} 
where 
$F_0$ 
is the diagonal operator defined by 
\begin{equation}
{F_0}{e_{m,n}} = 
\left\{
\begin{array}{cc}
{\frac{m+in}{ {(m^2 + n^2)}^{1/2}}}{e_{m,n}} &: (m,n) \neq 0 \cr
{e_{0,0}} &: (m,n) =(0,0).\cr
\end{array}
\right.
\end{equation}

 We note that the Baum-Connes assembly map [BCH94]
\begin{equation}
\mu : {KK^i}( C_0 ( B{\bf Z}^2 )),{\bf C}) = {KK^i}( C({\bf T}^2),{\bf C}) 
\rightarrow 
{K_i}( {C^{*}_r}({\bf Z}^2)) = {K_i}( C({\bf T}^2))
\end{equation}
 identifies the even Fredhom modules ${\bf z}_0$ and ${\bf Dirac}$ with $\pm [1]$ and $\pm [P_a]$ respectively, and the odd Fredholm modules 
 ${\bf w}_1$ and ${\bf w}_1 '$ with $\pm [U]$ and $\pm [V]$.
\end{proof}

\begin{Lemma}
\label{chernwonewoneprime}
The Chern characters of
${\bf w}_1$ and
${{\bf w}_1}'$
pair nontrivially with the generators
$[U]$ and $[V]$ of 
${K_1}(C( {\bf T}^2 ))$
as follows. We have 
 $<{ch}_{*}({\bf w}_1),[U]> = 1$, 
 $<{ch}_{*} ({\bf w}_1),[V]> = 0$, 
 $<{ch}_{*}({\bf w}_1 '),[U]> = 0$ and 
 $<{ch}_{*} ({\bf w}_1 '),[V]> = 1$.
\end{Lemma}
\begin{proof} 
 These calculations all follow from Example 2. 
\end{proof}

Via the *-homomorphism 
$\phi$,
 the Fredholm modules 
${\bf w}_1$ and 
${{\bf w}_1}'$
pull back to Fredholm modules 
\begin{equation}
{\bf z}_1 = {\phi^{*}}({\bf w}_1) = ( {l^2}({\bf Z}),{\pi_1} \circ \phi ,F),
\end{equation}
\begin{equation}
{{\bf z}_1}'= {\phi^{*}}({{\bf w}_1}') = ( {l^2}({\bf Z}),{\pi_1}' \circ \phi ,F)
\end{equation}
 over 
$C^{*} (H_3)$.
Explicitly, 
\begin{equation}
{\pi_1} \circ \phi (U) = S, \quad 
{\pi_1} \circ \phi (V) = I = {\pi_1} \circ \phi (W),
\end{equation}
\begin{equation}
{\pi_1}' \circ \phi (U) = I, \quad
{\pi_1}' \circ \phi (V) = S, \quad
{\pi_1}' \circ \phi (W) = I
\end{equation}
where $S$ is the shift operator. 

\begin{Prop}
\label{khomhbg1}
The Fredholm modules  ${\bf z}_1$ and  ${{\bf z}_1}'$  
are two linearly independent generators of  
 ${KK^1}( {C^{*}}({H_3}) , {\bf C})$.
\end{Prop}
\begin{proof} 
 It is immediate from the calculations of Lemma 6 that
\begin{equation}
< {ch}_{*}({\bf z}_1) , [U] > = 
< {ch}_{*}( {\phi^{*}}({\bf w}_1)) , [U] > =
< {ch}_{*}({\bf w}_1) , [ \phi(U)] > =
< {ch}_{*}({\bf w}_1) , [U] > = 1,
\end{equation}
\begin{equation}
< {ch}_{*} ({\bf z}_1) , [V] > = 0 
= < {ch}_{*} ({\bf z}_1) , [W] >,
\end{equation} 
\begin{equation}
 < {ch}_{*}({{\bf z}_1}') , [V] > = 1, \,\,\,
< {ch}_{*} ({{\bf z}_1}') , [U] > = 0
 = < {ch}_{*} ({{\bf z}_1}') , [W] >.
\end{equation}
To calculate the pairings of these Fredholm modules with the unitary
${V_{a}}$, we recall from (2) that  
\begin{equation}
V_a = {Z_a}( V \otimes {I_2}) {P_a} + (I_2 - P_a).
\end{equation}
It follows that 
\begin{equation}
 \phi(V_a) = \phi(Z_a) 
\left(
\begin{array}{cc}
V & 0 \cr
0 & 0 \cr
\end{array}
\right)
 +
\left(
\begin{array}{cc}
0 & 0 \cr
0 & 1 \cr
\end{array}
\right) 
=
 \left(
\begin{array}{cc}
V & 0 \cr
0 & 1 \cr
\end{array}
\right)
\end{equation}
since $\phi(Z_a) = I_2$. 
 Hence ${\phi_{*}} [V_a] =$ 
$[ \phi(V_a) ] = $
$[V] \in {K_1}(C ({\bf T}^2))$. 
 We have 
$\pi_1 (V) = I$, 
hence 
$(\pi_1 \otimes id)  \circ \phi (V_a) = I_2$ 
and it follows that 
\begin{equation*}
< {ch}_{*}({\bf z}_1) , [{V_{a}}] > =
 < {ch}_{*}({\phi^{*}}({\bf w}_1)) , [{V_{a}}] > =
< {ch}_{*}({\bf w}_1) , [\phi({V_{a}})] > =
\end{equation*}
\begin{equation}
< {ch}_{*}({\bf w}_1) , [I_2] > = 
2 < {ch}_{*}({\bf w}_1) , [1] > = 0.
\end{equation}
Since $\pi_1 ' (V) = S$, we have 
\begin{equation}
< {ch}_{*} ({{\bf z}_1}') , [{V_{a}}] > =
< {ch}_{*} ({{\bf z}_1}') , [V] > = 1.
\end{equation}
 Since all these pairings are either 0 or 1, it follows  from Connes' index theorem [Co94], p296, that 
${\bf z}_1$ and ${{\bf z}_1}'$ 
are generators, in the sense that if 
${\bf z}_1 = n {\bf w}$, 
for some $n \in {\bf Z}$ and some Fredholm module ${\bf w}$, 
then $n = \pm 1$. 
 This completes the proof of the proposition.
 Note that we will find a third generator of the odd K-homology in Corollary 17.
\end{proof}

We now describe even Fredholm modules.

\begin{Prop}
\label{khomhbg2}
The canonical even Fredholm module ${\bf z}_0$ is a nontrivial generator of\\ 
 ${KK^0}({C^{*}}(H_3), {\bf C})$, 
 which pairs to 1 with each of the generators 
 $[1]$, $[P_a]$ and $[P_b]$ of 
  ${K_0}( {C^{*}}(H_3))$. 
\end{Prop}
\begin{proof}
The *-homomorphism 
$\psi : {C^{*}}(H_3) \rightarrow {\bf C}$ 
defined on generators by
$\psi(U)=$
$\psi(V)=$
$\psi(W)=1$, 
 together with 
\begin{equation}
F=
\left(
\begin{array}{cc}
0 & 1 \cr
1 & 0 \cr
\end{array}
\right),
\gamma=
\left(
\begin{array}{cc}
1 & 0 \cr
0 & -1 \cr
\end{array}
\right)
\end{equation}
gives the canonical even Fredholm module
${\bf z}_0 = ({\bf C}^2, \psi \oplus 0 , F, \gamma)$ (Example 1).
It is immediate from Lemma 2 that 
$<{ch_{*}}( {\bf z}_0 ) , [1]> = 1$.
 Using the definition of the matrices $P_a$ and $P_b$ from [AP89]
 we calculate that 
\begin{equation}
\psi( P_a) = 
\left(
\begin{array}{cc}
1 & 0 \cr
0 & 0 \cr
\end{array}
\right)
= \psi( P_b)
\end{equation}
 hence 
$<{ch_{*}}( {\bf z}_0 ) , [{P_{a}}]> =$ 
$ <{ch_{*}}( {\bf z}_0 ), [1]> = 1$
$= <{ch_{*}}( {\bf z}_0 ) , [{P_{b}}]>$
as claimed.
\end{proof}

Consider again the quotient *-homomorphism $\phi : {C^{*}}(H_3) \rightarrow C({\bf T}^2)$ given by $U \mapsto U$, $V \mapsto V$, $W \mapsto I$.  
\begin{Lemma} 
\label{phikthy}
The induced map on K-theory
 ${\phi_{*}} : {K_0}({C^{*}}({H_3})) \rightarrow {K_0}(C({\bf T}^2))$ 
 satisfies
\begin{equation}
{\phi_{*}}[1] = [1] 
= {\phi_{*}}[{P_{a}}] 
= {\phi_{*}}[P_{b}].
\end{equation}
\end{Lemma}
\begin{proof}
 It is immediate that 
$ \phi(P_a) =
\left(
\begin{array}{cc}
1 & 0 \cr
0 & 0 \cr
\end{array}
\right)
 = \phi(P_b)$.
 Hence as elements of ${K_0}(C({\bf T}^2))$, we have 
$[\phi(P_a)] = [1] = [\phi(P_b)]$. 
And obviously $\phi(1) = 1$. 
\end{proof}

\begin{Lemma}
\label{khomhbg3}
Under the induced map on K-homology\\  
${\phi^{*}} : {KK^i}( C({\bf T}^2), {\bf C}) \rightarrow$
 ${KK^i}( {C^{*}}({H_3}), {\bf C})$ 
 we have 
${\phi^{*}}({\bf w}_0 ) ={\bf z}_0$. 
We also obtain a Fredholm module 
${\phi^{*}}( {\bf Dirac})$ 
which pairs to zero with ${K_0}( {C^{*}}({H_3}) )$. 
\end{Lemma}
\begin{proof} The fact that ${\phi^{*}}({\bf w}_0 ) ={\bf z}_0$ follows immediately from the definition of these Fredholm modules.  
 We describe the Fredholm module 
${\phi^{*}}( {\bf Dirac})$. 
 Consider ${C^{*}}(H_3)$ acting on 
${l^2}({\bf Z}^2)$ via  
\begin{equation}
U{e_{m,n}} = e_{m+1,n}, \quad
V{e_{m,n}} = e_{m,n+1}, \quad
W{e_{m,n}} = {e_{m,n}}.
\end{equation}  
Take $H = {l^2}({\bf Z}^2) \oplus {l^2}({\bf Z}^2)$, 
and $\pi \circ \phi : {C^{*}}(H_3) \rightarrow {\bf B}(H)$ to be 
\begin{equation}
\pi \circ \phi (a) = \left(
\begin{array}{cc}
a & 0 \cr
0 & a \cr
\end{array}
\right).
\end{equation}
 Finally, take $F$ as defined in (22). 
Then 
$\phi^{*}( {\bf Dirac} ) = \left(
 H, \pi \circ \phi, F 
\right)$.
 Using Lemma 9 we have 
\begin{equation}
<{ch_{*}}({\phi^{*}}({\bf Dirac})),[1]> 
= <{ch_{*}}({\bf Dirac}),{\phi_{*}}[1]> 
=  <{ch_{*}}({\bf Dirac}),[1]> 
= 0,
\end{equation}
since to calculate this we need to consider expressions of the form 
$Tr( \gamma \pi(1) {[F, \pi(1)]}^{2k} )$, which are all zero since $[F, \pi(1)] = 0$. 
By the same argument
\begin{equation}
<{ch_{*}}({\phi^{*}}({\bf Dirac})),[{P_{a}}]> = 0 = 
<{ch_{*}}({\phi^{*}}({\bf Dirac})),[{P_{b}}]>.
\end{equation}
 We would like to deduce from this that 
${\phi^{*}}({\bf Dirac}) = {\bf 0}$. However, at this point we do not  a priori know whether the pairing between K-theory and K-homology given by 
$<.,.>$ together with the Chern character is faithful, although we will see later that this is true.  
\end{proof}

\section{Six term cyclic sequence for K-homology}

We now consider the six term cyclic exact sequence for K-homology of crossed products by ${\bf Z}$. 
 Associated to any crossed product algebra $A \times_{\alpha} {\bf Z}$  
is the following semisplit short exact sequence of C*-algebras, 
the Pimsner-Voiculescu ``Toeplitz extension" 
[PV80]:
\begin{equation}
\label{toeplitzext}
 0 \rightarrow {A \otimes {\bf K}} \rightarrow {T_{\alpha}} \rightarrow  {A \times_{\alpha} {\bf Z}} \rightarrow  0.
\end{equation}
 Here ${T_{\alpha}}$ 
is the C*-subalgebra of 
$( {A {\times}_{\alpha} {\bf Z}}) \otimes {\it T} $
generated by 
$a \otimes 1$, 
$a \in A$ 
and 
$ V \otimes f$,
where 
$V$ 
is the unitary implementing the action of 
$\alpha$ 
on 
$A$, 
and 
$f$ 
is the non-unitary isometry generating the ordinary Toeplitz algebra
 $T$, 
that is  
$f \in {\bf B}({l^2}({\bf N}))$,
$f{e_n} = {e_{n+1}}$.
 This extension defines the Toeplitz element 
${\bf x} \in {KK^1}( {A \times_{\alpha} {\bf Z}}, A) $.  

Applying the K-functor gives the Pimsner-Voiculescu six term cyclic sequence for K-theory. 
The corresponding six term cyclic sequence for K-homology is:

\begin{diagram}
{{KK^0}( A  , {\bf C})} & \lTo^{id - {\alpha}^{*}} & {{KK^0}( A, {\bf C})} & \lTo^{i^{*}} & {{KK^0}(A {\times_{\alpha}} {\bf Z}, {\bf C})}  \\
\dTo^{\partial_0} & & & & \uTo^{\partial_1} \\
{{KK^1}( A {\times_{\alpha}} {\bf Z} , {\bf C})} & \rTo^{i^{*}} & {{KK^1}( A, {\bf C})} & \rTo^{id - {\alpha}^{*}} & {{KK^1}( A, {\bf C} )}\\
\end{diagram}

Here $i$ denotes the inclusion map 
 $i : A \hookrightarrow A \times_{\alpha} {\bf Z}$. 
The vertical maps 
${\partial_0}$ and  ${\partial_1}$  
are given by taking the Kasparov product with the element 
${\bf x} \in {KK^1}( A \times_{\alpha}{\bf Z}, A)$ 
corresponding to the Toeplitz extension (43).
 This sequence formulated in terms of Ext appears in the original paper of Pimsner and Voiculescu [PV80].
However, the relationship between Ext and the Fredholm module picture of K-homology is not transparent. 

Let $A$ be a unital C*-algebra, $\phi : A \rightarrow {\bf C}$ a nonzero  *-homomorphism, and  ${\bf z}_0$ and ${\bf z}_1$ the canonical even and odd Fredholm modules associated to $\phi$ defined in Examples 1 and 2.
We have:

\begin{Prop}
\label{bdy0}
[Ha01]
 $\partial_0 ( {\bf z}_0 ) = {\bf z}_1$.
\end{Prop}

For the discrete Heisenberg group, we take 
$A= C( {\bf T}^2 )$, generated by commuting unitaries $U$ and $W$, together with a ${\bf Z}$-action implemented by an automorphism $\alpha$ given by $\alpha(U) = WU$, $\alpha(W)=W$.
 Then as we saw previously $A \times_{\alpha} {\bf Z} \cong {C^{*}}(H_3)$.  
 Note that we are changing our notation for $C({\bf T}^2)$ slightly from the last section, since now we are considering it to be a subalgebra, rather than a quotient, of ${C^{*}}(H_3)$. 

\begin{Prop}
\label{khomhbgexact1}
Under the map
 $(id - {\alpha}^{*}) : {KK^1}( C( {\bf T}^2 ), {\bf C} ) \rightarrow {KK^1}( C( {\bf T}^2 ), {\bf C} )$
 we have 
$(id - {\alpha^{*}}) ({\bf w}_1) = {\bf 0} $ and  
$(id - {\alpha^{*}}) ({{\bf w}_1}') = - {\bf w}_1$.
\end{Prop}
\begin{proof}   
 We saw in Prop 5 that the odd K-homology   
${KK^1}( C( {\bf T}^2 ),{\bf C}) \cong {\bf Z}^2$ 
is generated by the two Fredholm modules 
${\bf w}_1 = ({l^2}( { \bf Z} ), {\pi_1}, F)$, 
${{\bf w}_1}' = ({l^2}( { \bf Z} ), {\pi_1}', F)$,
where the representations 
${\pi_1},$ 
 ${\pi_1}'$ 
of 
$C( {\bf T}^2 )$ 
on 
${l^2}( { \bf Z} )$
 are given by
\begin{equation}
{\pi_1}(U) = S, \quad
{\pi_1}(W) = I, \quad
{\pi_1}'(U) = I, \quad
{\pi_1}'(W) = S.
\end{equation}
Here 
$S$ 
is the shift operator with respect to the orthonormal basis
$ { \{ {e_n} \} }_{ n \in {\bf Z}}$,
$S{e_n} = {e_{n+1}}$,
and 
$F{e_n}= sign(n){e_n}$.
Then 
\begin{equation}
{\pi_1}( \alpha(U) ) 
= {\pi_1}(UW) 
= {\pi_1}(U),\quad
 {\pi_1}( \alpha(W)) 
= {\pi_1}(W).
\end{equation}
Hence 
${\alpha^{*}} ({\bf w}_1 ) = {\bf w}_1$ as Fredholm modules,
so 
$ (id - {\alpha}^{*}) ( {\bf w}_1) = {\bf 0}$.
Similarly, 
\begin{equation}
{\pi_1}'( \alpha(U) ) 
= {\pi_1}'( WU ) 
= {\pi_1}'(W),\quad
{\pi_1}'( \alpha(W) ) 
= {\pi_1}'(W).
\end{equation}
 So 
${\alpha^{*}}({{\bf w}_1}') = {{\bf w}_1} + {{\bf w}_1}' $,
and hence  
$(id - {\alpha}^{*})({{\bf w}_1}') = - {{\bf w}_1}$.  
So  
$Ker(id - {\alpha}^{*}) \cong {\bf Z}$, generated by ${\bf w}_1$, and 
 $Im(id - {\alpha}^{*}) \cong {\bf Z}$, generated by ${\bf w}_1$.
\end{proof}

\begin{Prop}
\label{idminusalphastar}
The map
 $ (id - {\alpha}^{*}) : {KK^0}( C( {\bf T}^2 ), {\bf C} ) \rightarrow {KK^0}( C( {\bf T}^2 ), {\bf C} )$
 is the zero map.
\end{Prop}  
\begin{proof} 
We saw in Prop 5 that the even K-homology
${KK^0}( C({\bf T}^2),{\bf C}) \cong {\bf Z}^2 $
is generated by the 
canonical even Fredholm module 
${\bf w}_0 = ({\bf C}^2, \phi \oplus 0 , F, \gamma)$
corresponding to the *-homomorphism 
 $\phi : C({\bf T}^2) \rightarrow {\bf C}$ given by 
$\phi(U) = \phi (W) =1$, 
 together with the Fredholm module
${\bf Dirac}$, 
corresponding to the phase of the Dirac operator on
${\bf T}^2$.
 Now, ${\alpha^{*}}( {\bf w}_0) = ({\bf C}^2, (\phi \circ \alpha) \oplus 0, F, \gamma)$, 
and we check on generators that $\phi \circ \alpha (U) = \phi(WU) = 1 = \phi(U)$, 
and 
$\phi \circ \alpha (W) = \phi (W)$.
Hence $\phi \circ \alpha = \phi$, and so 
${\alpha^{*}}( {\bf w}_0) = {\bf w}_0$ as Fredholm modules. 

\begin{Lemma} ${\alpha^{*}}( {\bf Dirac}) ={\bf Dirac}$ as elements of ${KK^0}( C({\bf T}^2), {\bf C})$. 
\end{Lemma}
\begin{proof} 
We have 
\begin{equation}
{\alpha^{*}}( {\bf Dirac}) = ( {l^2}({\bf Z}^2) \oplus {l^2}({\bf Z}^2), ({\pi_0}  \oplus {\pi_0}) \circ \alpha, 
F = \left(
\begin{array}{cc}
0 & {F_0}\cr
{F_0} & 0 \cr
\end{array}
\right)
)
\end{equation} 
 where 
\begin{equation}
{\pi_0} \circ \alpha(U) {e_{m,n}} = {\pi_0}(UW) e_{m,n} = e_{m+1,n+1},\quad
{\pi_0} \circ \alpha(W) e_{m,n} = {\pi_0}(W) e_{m,n} = e_{m,n+1}, 
\end{equation}
 and 
\begin{equation}
{F_0}{e_{p,q}} = 
\left\{
\begin{array}{cc}
{\frac{p+iq}{(p^2 + q^2)^{1/2}}}{e_{p,q}} &: (p,q) \neq (0,0) \cr 
 e_{0,0} &: (p,q)=(0,0). \cr
\end{array}
\right. 
\end{equation}
 Recall that Fredholm modules $(H_0, \pi_0, F_0, \gamma_0)$ and 
$(H_1, \pi_1, F_1, \gamma_1)$ over a C*-algebra $A$ 
are said to be unitarily equivalent 
if there exists a unitary $T : H_0 \rightarrow H_1$ such that 
${T^{*}} \pi_1 (a) T = \pi_0 (a)$, for each $a \in A$, 
${T^{*}} F_1 T = F_0$, and ${T^{*}} \gamma_1 T = \gamma_0$. 
Unitarily equivalent Fredholm modules define the same element of K-homology. 
 
Consider the unitary ${T_0} \in {\bf B}( {l^2}({\bf Z}^2))$ defined by 
${T_0} e_{m,n} = e_{m,n-m}$. 
It is easy to check that ${T_0}^{*} {\pi_0}(U) {T_0} = {\pi_0}(UW)$ and ${T_0}^{*} {\pi_0}(W) T_0 = {\pi_0}(W)$. 
 Taking 
\begin{equation}
T = 
\left(
\begin{array}{cc}
{T_0} & 0 \cr
0 & {T_0} \cr
\end{array}
\right)
\in {\bf B} (H),
\end{equation}
 we see that ${\alpha^{*}}({\bf Dirac})$ is unitarily equivalent via $T$ 
to a  Fredholm module 
\begin{equation}
{\bf y}_1 ' = 
(
H, {\pi_1} \oplus {\pi_1}, \tilde{F_1} = 
\left(
\begin{array}{cc}
0 & {F_1} \cr
{F_1}^{*} & 0 \cr
\end{array}
\right)
) 
\end{equation}
 where ${\pi_1}(U) e_{p,q} = e_{p+1,q}$, 
${{\pi_1}(W)} e_{p,q} = e_{p,q+1}$, and
  $F_1$ is defined by :
\begin{equation}
{F_1} e_{p,q} = 
\left\{
\begin{array}{cc}
{\frac{p + i(p +q)}{(p^2 + (p+q)^2)^{1/2}}} e_{p,q} & : (p,q) \neq (0,0) \cr
e_{0,0} & : (p,q)=(0,0) \cr
\end{array}
\right.
\end{equation}
 We construct a homotopy of Fredholm modules from ${\bf y}_1 '$ to ${\bf Dirac}$. 
For $0 \leq t \leq 1$, we define 
\begin{equation}
{\bf y}_t ' = 
(
H, {\pi_1} \oplus {\pi_1}, \tilde{F_t} = 
\left(
\begin{array}{cc}
0 & {F_t} \cr
{F_t}^{*} & 0 \cr
\end{array}
\right)
) 
\end{equation}
 where 
$F_t$ is given by :
\begin{equation}
{F_t} e_{p,q} = 
\left\{
\begin{array}{cc}
{\frac{p + i(tp +q)}{(p^2 + (tp+q)^2 )^{1/2}}} e_{p,q} & : (p,q) \neq (0,0) \cr
e_{0,0} & : (p,q)=(0,0) \cr
\end{array}
\right.
\end{equation}
 So 
${ \{ F_t \} }_{0 \leq t \leq 1}$ 
is a norm continuous path of elements of 
${\bf B}(H)$, 
and
  ${\bf y}_0 ' ={\bf Dirac}$. 
Hence 
 ${\bf Dirac}$ and ${\alpha^{*}}({\bf Dirac})$  define the same element of 
 ${KK^0}(C({\bf T}^2),{\bf C})$. 
\end{proof}

So $(id - \alpha^{*})$ is the zero map as claimed, thus completing the proof of  Prop 13.
\end{proof}

\begin{Cor}
\label{partial0cor}
The image of the map 
$\partial_0  :  {KK^0}(C({\bf T}^2), {\bf C}) \rightarrow$
$ {KK^1}( {C^{*}}({H_3}), {\bf C})$, 
 ${\bf z} \mapsto {\bf x}{\hat{\otimes}} {\bf z}$
 is a copy of ${\bf Z}^2$.
\end{Cor}

 \begin{Prop} 
\label{partial0prop} 
 $\partial_0 ({\bf w}_0) = {{\bf z}_1}'$, ${i^{*}}( {{\bf z}_1}')= {\bf 0}$, 
${i^{*}}( {{\bf z}_1})= {\bf w}_1$. 
\end{Prop} 
 \begin{proof}  
Applying Prop 11, we see that
$\partial_0 ({\bf w}_0)$ 
 is the Fredholm module 
 $({l^2}({\bf Z}) , \pi, F) \in$\\
 ${KK^1}( {C^{*}}({H_3}),{\bf C})$, 
where $ \pi : {C^{*}}( H_3 ) \rightarrow {\bf B}( {l^2}({\bf Z}))$ is given by 
\begin{equation}
(\pi(U) \xi) (n) = \psi( {\alpha^n} (U) ) \xi (n) = \psi( {W^n}U ) \xi (n) =
\xi(n),
\end{equation}
\begin{equation}
(\pi(W) \xi) (n) = \psi( {\alpha^n} (W) ) \xi (n) = \psi(W ) \xi (n) =
\xi(n),
\end{equation}
\begin{equation}
(\pi(V) \xi) (n) = \xi(n+1),
\end{equation}
 for $\xi \in {l^2}({\bf Z})$, and we have $F \xi(n) = sign(n) \xi(n)$.
 Hence $\partial_0 ({\bf w}_0) = {{\bf z}_1}'$ as Fredholm modules. 

We have  
${i^{*}}( {{\bf z}_1}') = 
( {l^2}({\bf Z}),  {\pi_1 '} \circ \phi \circ i , F)
 \in {KK^1}( C({\bf T}^2), {\bf C})$.
 Since ${\pi_1 '}\circ \phi \circ i (U) = I = {\pi_1 '} \circ \phi \circ i (W)$ in ${\bf B}( {l^2}({\bf Z}))$, it follows that 
 ${i^{*}}( {{\bf z}_1}')$ is a degenerate Fredholm module, and hence is a trivial element of K-homology. 
 Further, 
 we have  that
${i^{*}}( {{\bf z}_1}) = ( {l^2}({\bf Z}),  \pi_1 \circ \phi \circ i , F)$.
 Now, 
${\pi_1}\circ \phi \circ i (U) = S$ and 
${\pi_1}\circ \phi \circ i (W) = I$.
 So ${i^{*}}( {{\bf z}_1}) = {\bf w}_1$ as Fredholm modules. 
 It follows that 
${{\bf z}_1}$
is not in the image of 
$\partial_0$.
\end{proof}

\begin{Cor}
\label{thirdgen}
 A third element of ${KK^1}( {C^{*}}(H_3), {\bf C})$, not in the span of ${\bf z}_1$ and ${\bf z}_1 '$, is given by the Fredholm module ${\partial_0}( {\bf Dirac})$, which we do not describe explicitly. 
\end{Cor}
\begin{proof}
 Two generators are given by ${\bf z}_1$ and ${\bf z}_1 '$ (Prop 7) and it follows from Corollary 15 and Prop 16 that 
${\partial_0}( {\bf Dirac})$ 
is a nontrivial element of K-homology not in the span of 
${\bf z}_1$ and ${\bf z}_1 '$.
Hence the result. 
\end{proof}

The pairings of ${\partial_0}( {\bf Dirac})$ with K-theory are as follows: 

\begin{Lemma}
$< {ch_{*}}( {\partial_0}( {\bf Dirac})), [U] > = 0 = < {ch_{*}}( {\partial_0}( {\bf Dirac})), [V] >$, and\\ 
 $< {ch_{*}}( {\partial_0}( {\bf Dirac})), [V_a] > =1$.
\end{Lemma}
\begin{proof}
 We use the duality of $\partial_0$ with the Pimsner-Voiculescu boundary map\\ 
  $\delta_1 : {K_1}( {C^{*}}(H_3) ) \rightarrow {K_0}( C( {\bf T}^2))$, 
 a map which satisfies [AP89] 
$ \delta_1 [U] = 0$, 
 $\delta_1 [V] = [1]$ and 
 $\delta_1 [V_a] = [P_a]$. 
 It follows that 
\begin{equation}
< {ch_{*}}( {\partial_0}( {\bf Dirac})), [U] > = < {ch_{*}}( {\bf Dirac}), 0 > = 0,
\end{equation}
\begin{equation}
< {ch_{*}}( {\partial_0}( {\bf Dirac})), [V] > = < {ch_{*}}( {\bf Dirac}), [1] > = 0,
\end{equation}
\begin{equation}
< {ch_{*}}( {\partial_0}( {\bf Dirac})), [V_a] > = < {ch_{*}}( {\bf Dirac}), [P_a] > = 1.
\end{equation}
 Since these pairings are all either 0 or 1 it follows that ${\partial_0}( {\bf Dirac})$ is a generator of \\
 ${KK^1}( {C^{*}}( H_3), {\bf C})$. 
\end{proof}

We summarize our knowledge of the odd K-homology in :

\begin{Thm}
\label{oddKhomology}
The Fredholm modules 
${\bf z}_1$, ${\bf z}_1 '$ 
and 
${\partial_0}( {\bf Dirac})$ generate 
  ${KK^1}( {C^{*}}(H_3), {\bf C}) \cong {\bf Z}^3$. 
\end{Thm}

We give further verification of these results
by exploiting the duality of the maps 
$\partial_0$ and  $\partial_1$ 
with the Pimsner-Voiculescu boundary maps 
$\delta_i$
on K-theory. 
Given 
${\bf q} \in {KK^1}({\bf C}, {C^{*}}({H_3}) )$,
${\bf z} \in {KK^0}(C({\bf T}^2), {\bf C})$
 and the element ${\bf x} \in {KK^1}( {C^{*}}(H_3), C({\bf T}^2))$ corresponding to the Toeplitz extension (43),   
by associativity of the Kasparov product we have
\begin{equation}
{\bf q} {\hat{\otimes}}( {\bf x} {\hat{\otimes}} {\bf z}) =
( {\bf q} {\hat{\otimes}} {\bf x}) {\hat{\otimes}} {\bf z}.
\end{equation}
So
\begin{equation} 
< {ch_{*}}( \partial_0 ({\bf z})), {\bf q} >
 = < {ch_{*}}( {\bf z}), \delta_1 ({\bf q}) >.
\end{equation}
The K-group 
${K_1}({C^{*}}({H_3})) \cong $  
${KK^1}( {\bf C}, {C^{*}}({H_3}))$ 
is generated by 
$[U]$, 
$[V]$ 
and  
$[{V_{a}}]$, 
and we know from (6) that under the 
Pimsner-Voiculescu boundary map
\begin{equation}
\delta_1  :  {KK^1}( {\bf C}, {C^{*}}({H_3})) \rightarrow {KK^0}( {\bf C}, C({\bf T}^2)), 
 {\bf q} \mapsto {\bf q } {\hat{\otimes}} {\bf x}
\end{equation}
 we have  
${\delta_1}[U] = 0$,
${\delta_1}[V] = [1]$ 
and 
${\delta_1}[{V_{a}}]=[{P_{a}}]$, 
the Bott generator for 
${K_0}(C({\bf T}^2))$. 
So for any 
${\bf w} \in {KK^0}( C({\bf T}^2), {\bf C})$,
\begin{equation}
< {ch_{*}}( \partial_0 ({\bf w})), [U] > 
= [U] {\hat{\otimes}} ( {\bf x} {\hat{\otimes}} {\bf w}) 
= {\delta_1}[U] {\hat{\otimes}} {\bf w} = 0
\end{equation}
while we know that
$ < {ch_{*}}({{\bf z}_1}), [U] > = 1$.
Hence
$ {{\bf z}_1} \neq \partial_0 ({\bf w})$, for any ${\bf w} \in$
 ${KK^0}( C({\bf T}^2), {\bf C})$. 
We check further that 
\begin{equation}
1 = < {ch_{*}}({{\bf z}_1}' ), [{V_{a}}] > 
= [{V_{a}}] {\hat{\otimes}} {\partial_0 ( {\bf w}_0 )} 
= \delta_1 ( {V_{a}} ) {\hat{\otimes}} {{\bf w}_0} 
= < {ch_{*}}({{\bf w}_0}), [{P_{a}}] > = 1,
\end{equation}
\begin{equation}
1 = < {ch_{*}}({{\bf z}_1}' ), [V] > 
= [V] {\hat{\otimes}} {\partial_0 ( {\bf w}_0 )} 
= \delta_1 (V ) {\hat{\otimes}} {{\bf w}_0} 
= < {ch_{*}}({{\bf w}_0}), [1] > = 1,
\end{equation}
\begin{equation}
0 = < {ch_{*}}({{\bf z}_1}' ), [U] > 
= < {ch_{*}}({{\bf w}_0}), \delta_1 [U] > 
= < {ch_{*}}({{\bf w}_0}), {\bf 0} > = 0.
\end{equation}

\begin{Prop}
\label{khomhbgexact2} 
 The image $\partial_1 ( {\bf w}_1 )$ of the Fredholm module ${\bf w}_1$ under the map\\ 
 $\partial_1  :  {KK^1}(C({\bf T}^2), {\bf C}) \rightarrow {KK^0}( {C^{*}}({H_3}), {\bf C})$
 pairs identically to zero with ${K_0}( {C^{*}}({H_3}))$. 
\end{Prop}
\begin{proof} 
The image 
$\partial_1 ({\bf w}_1)$ 
of 
${\bf w}_1$ is the even Fredholm module  
\begin{equation}
\partial_1 ({\bf w}_1) = 
({l^2}({\bf Z}^2) {\oplus} {l^2}({\bf Z}^2), {\pi \oplus \pi}, 
F = 
\left(
\begin{array}{cc}
0 & {F_0} \cr
{F_0}^{*} & 0 \cr
\end{array}
\right)
)
\end{equation}
with
$\pi(U){e_{m,n}} = {e_{m+1,n}}$, 
$\pi(V){e_{m,n}} = {e_{m,n+1}}$ and 
$\pi(W) = I$. 
We define 
${F_0}$ 
by
\begin{equation}
{F_0}{e_{m,n}} = 
\left\{
\begin{array}{cc}
{\frac{m+in}{{(m^2 + n^2)}^{1/2}}}{e_{m,n}} & : (m,n) \neq (0,0) \cr
e_{0,0} & : (m,n) = (0,0) \cr
\end{array}
\right.
\end{equation} 
 We calculate the pairings of $\partial_1 ({\bf w}_1)$ with the generators 
$[1]$, $[P_a]$ and $[P_b]$ of 
${K_0}( {C^{*}}(H_3))$ 
by exploiting the duality of $\partial_1$ with the Pimsner-Voiculescu boundary map $\delta_0$ on K-theory 
\begin{equation}
\delta_0 : {K_0}( C^{*} (H_3) ) \rightarrow {K_1}( C( {\bf T}^2))
\end{equation}
 We have from (4) that   
 $\delta_0 [1] = 0$, 
 $\delta_0 [P_a] = 0$ and 
 $\delta_0 [P_b] = [W]$. 
 It is immediate that\\   
 $<{ch_{*}}(\partial_1 ({\bf w}_1)),[1]> = 0$, 
 since $[F, ( \pi \oplus \pi ) (1) ] = 0$. 
 To calculate  
$<{ch_{*}}(\partial_1 ({\bf w}_1)),[P_a]>$, we note that 
$\pi( P_a ) = 
\left(
\begin{array}{cc}
1 & 0 \cr
0 & 0 \cr
\end{array}
\right)$. 
 Hence 
\begin{equation}
<{ch_{*}}(\partial_1 ({\bf w}_1)),[P_a]> = <{ch_{*}}(\partial_1 ({\bf w}_1)),[1]> = 0.
\end{equation}
 We also have    
\begin{equation}
<{ch_{*}}(\partial_1 ({\bf w}_1)),[P_b]> = 
<{ch_{*}}({\bf w}_1), \delta_0 [P_b]> = 
<{ch_{*}}({\bf w}_1), [W]> = 0. 
\end{equation}
 We would like to deduce from this 
 that 
$\partial_1 ({\bf w}_1) = {\bf 0}$, but again we do not know that $<.,.>$ is faithful.
\end{proof}

\begin{Lemma}
The Fredholm module $\partial_1 ( {\bf w}_1 ')$ is a nontrivial element of 
 ${KK^0}( {C^{*}}(H_3), {\bf C})$ 
not in the span of ${\bf z}_0$. 
\end{Lemma}
 \begin{proof}
Again we use the duality with the Pimsner-Voiculescu map.   
 We have 
\begin{equation}
<{ch_{*}}(\partial_1 ({{\bf w}_1}')),[1]> = 
 <{ch_{*}}({{\bf w}_1}'),\delta_0 [1]> = 0.
\end{equation}
 It follows that 
$<{ch_{*}}(\partial_1 ({{\bf w}_1}')),[{P_{a}}]> = 0$, and furthermore 
 \begin{equation} 
  <{ch_{*}}(\partial_1 ({{\bf w}_1}')),[{P_{b}}]> = 
<{ch_{*}}({{\bf w}_1}'),[W]> = 1.
\end{equation}
 Hence $\partial_1 ({{\bf w}_1}')$ is a nontrivial element of K-homology. 
 Since by Prop 8 we know that the canonical even Fredholm module ${\bf z}_0$ pairs to 1 with each of $[1]$, $[P_a]$ and $[P_b]$, it follows that 
  $\partial_1 ({{\bf w}_1}')$ and ${\bf z}_0$ are linearly independent.
\end{proof}

We describe a third even Fredhom module. 
It follows from  Prop 13 that the map
\begin{equation}
{i^{*}} : {KK^0}( {C^{*}}(H_3) , {\bf C}) \rightarrow {KK^0}( C({\bf T}^2) , {\bf C})
\end{equation}
 is surjective. 
It is immediate that  ${i^{*}}({\bf z}_0) = {\bf w}_0$. 
 Denote by ${\bf Dirac'}$ an element of\\
  ${KK^0}( C^{*}(H_3) , {\bf C})$ 
satisfying ${i^{*}}({\bf Dirac'}) = {\bf Dirac}$. 
We do not give an explicit description of such a Fredholm module, but we can calculate its pairings with the generators of K-theory. 

\begin{Lemma} 
$< {ch_{*}}( {\bf Dirac'}), [1] > =0$,  
$< {ch_{*}}( {\bf Dirac'}), [P_a] > =1$.
\end{Lemma}
\begin{proof}
We know that the map ${i_{*}} : {K_0}(C({\bf T}^2)) \rightarrow {K_0}( {C^{*}}(H_3))$ is injective. 
We use the same notation $[1]$ and $[P_a]$ for the generators of ${K_0}(C({\bf T}^2))$ and the corresponding generators of ${K_0}( {C^{*}}(H_3))$.  
We have 
\begin{equation*}
< {ch_{*}}( {\bf Dirac'}), [1] > = < {ch_{*}}( {\bf Dirac'}), {i_{*}}[1] >
= < {ch_{*}}( {i^{*}}({\bf Dirac}')), [1] >
\end{equation*}
\begin{equation} 
= < {ch_{*}}( {\bf Dirac}), [1] > = 0,
\end{equation}
\begin{equation}
< {ch_{*}}( {\bf Dirac'}), [P_a] > = < {ch_{*}}( {i^{*}}({\bf Dirac}')), [P_a] > =
< {ch_{*}}( {\bf Dirac}), [P_a] > = 1.
\end{equation}
Note that 
$< {ch_{*}}( {\bf Dirac'}), [P_b] >$
 is not uniquely defined, 
 because for any $k \in {\bf Z}$, 
we have \\
 ${i^{*}}( {\bf Dirac}' + k {\partial_1}( {\bf w}_1 ')) =$
 ${\bf Dirac}$, 
and  
$< {ch_{*}}({\partial_1}( {\bf w}_1 ')  ), [P_b] >=1$.
 However we do know that by Connes' index theorem [Co94], p296,  
$< {ch_{*}}( {\bf Dirac'}), [P_b] > \in {\bf Z}$. 
Recall that we previously found that the pairings of the Fredholm module
${\partial_1}( {\bf w}_1 ')$ with the generators 
$[1]$, $[P_a]$ and $[P_b]$ of K-theory are 0, 0 and 1 respectively. 
 Hence we can modify ${\bf Dirac'}$ by adding on multiples of 
${\partial_1}( {\bf w}_1 ')$ to ensure that 
$< {ch_{*}}( {\bf Dirac'}), [P_b] > = 0$. 
Hence ${\bf Dirac'}$ is a third generator of the even K-homology.
\end{proof}

We summarize our knowledge of the even K-homology in the following :

\begin{Thm}
\label{evenKhomology}
The Fredholm modules ${\bf z}_0$, ${\bf Dirac}'$ and  ${{\bf w}_1}'$ generate
  ${KK^0}( {C^{*}}( H_3), {\bf C}) \cong$  ${\bf Z}^3$.
\end{Thm}  
\begin{Prop}
\label{summaryhbgkthykhompairings}
We collect all our results on the pairings between the K-homology and K-theory of ${C^{*}}(H_3)$ in the following table.

\begin{equation*}
\begin{array}{cccc}
Even\,\, case  & {\bf z}_0 & {\bf Dirac'} & \partial_1( {\bf w}_1 ') \cr 
[1]  & 1 & 0 & 0 \cr
[P_a]  & 1 & 1 & 0 \cr
[P_b]  & 1 & 0 & 1 \cr
\end{array}
\end{equation*}
\begin{equation*}
\begin{array}{cccc}
Odd\,\, case & {\bf z}_1 & {\bf z}_1 ' & {\partial_0}( {\bf Dirac})  \cr 
[U] & 1 & 0 & 0 \cr
[V] & 0 & 1 & 0 \cr
[V_a] & 0 & 1 & 1 \cr
\end{array}
\end{equation*}
\end{Prop}

We see from this table that for $C^{*} (H_3)$ the pairing between K-theory and K-homology via the Chern character and $<.,.>$ is faithful. 
Given a Fredholm module ${\bf z} \in {KK^i}( C^{*} (H_3), {\bf C})$, if for every ${\bf q} \in$ $K_i ( C^{*} (H_3))$ we have 
$< {ch_{*}}({\bf z}), {\bf q}>=0$, then ${\bf z}$ is a trivial element of K-homology.

In conclusion, we have been able to give a  detailed description of the even and odd K-homology of ${C^{*}}(H_3)$. 
The six term cyclic sequence on K-homology played a central role in allowing us to construct Fredholm modules over ${C^{*}}(H_3)$ from known ones over $C({\bf T}^2)$.  
 The techniques we use generalize very nicely to other crossed products by ${\bf Z}$. 
This is the subject of [Ha01]
 where we study the morphisms
\begin{equation}
 {\partial_i} : {KK^i}(A, {\bf C}) \rightarrow {KK^{i+1}}( A \times_{\alpha} {\bf Z}, {\bf C})
\end{equation} 
 that  played a crucial role in this section. 
Conversely, for the infinite dihedral group, the methods we use have so far resisted generalization. 
 Given an algebra $A$ with a ${\bf Z}_2$-action $\sigma$, there is no apparent natural map from the K-homology of $A$ to the K-homology of the crossed product $A \times_{\sigma} {\bf Z}_2$, even for very restricted classes of algebras. 

As a footnote to the above, we compare the K-homology of ${C^{*}}(H_3)$ with that of ${C^{*}}( {\bf Z}^3 )$. 
It follows from the universal coefficient theorem (Prop 4) that:
\begin{equation}
 {KK^0}( {C^{*}}( {\bf Z}^3) , {\bf C}) \cong {\bf Z}^4, \quad
{KK^1}( {C^{*}}( {\bf Z}^3) , {\bf C}) \cong {\bf Z}^4 . 
\end{equation}
 For the generators of the even K-homology, we first of all have a canonical Fredholm module ${\bf z}_0$ (Example 1) corresponding to the *-homomorphism 
${C^{*}}( {\bf Z}^3) \rightarrow {\bf C}$ given by $U_1$, $U_2$, $U_3$ $ \mapsto 1$.
 We also have *-homomorphisms
$\phi_1$, $\phi_2$, 
$\phi_3 : {C^{*}}( {\bf Z}^3) \rightarrow {C^{*}}( {\bf Z}^2)$
 given by 
$\phi_i (U_i) = 1$, $\phi_i (U_j) = U_j$ for $j \neq i$. 
 Let ${\bf Dirac} \in {KK^0}( {C^{*}}( {\bf Z}^2), {\bf C})$ be the Fredholm module that we defined previously in (20). 
 Then the remaining three generators of ${KK^0}( {C^{*}}( {\bf Z}^3), {\bf C})$
are given by ${\phi_i^{*}}( {\bf Dirac})$, $i=1,2,3$. 

For the odd K-homology, recall that the odd K-homology of ${C^{*}}({\bf Z}) \cong$ $C({\bf T})$ is generated by a Fredholm module 
\begin{equation}
{\bf w}_1 = ( {l^2}({\bf Z}), \pi, F)
\end{equation}
 with 
$\pi (U) e_n = e_{n+1}$, 
and $F e_n = sign (n) e_n$. 
 We have *-homomorphisms 
$\psi_i : {C^{*}}( {\bf Z}^3 ) \rightarrow {C^{*}}( {\bf Z} )$, 
$i=1,2,3$ given by
$\psi_i (U_i) = U$, $\psi_i (U_j) =1$, for $i \neq j$. 
 Then three of the generators of 
${KK^1}( {C^{*}}( {\bf Z}^3) , {\bf C})$
 are given by 
${\psi_i^{*}}( {\bf w}_1 )$, $i=1,2,3$. 
 The fourth generator is given by ${\bf Dirac}$, the Fredholm module given by the bounded formulation of the Dirac operator on ${\bf T}^3$. 

In this example, the K-homology and the K-theory of ${C^{*}}( {\bf Z}^3)$ can be very nicely identified, using the assembly map. 
 Recall [BCH94] that for a general discrete group $\Gamma$, we have 
\begin{equation}
\mu : {KK^i}( C_0 ( B \Gamma), {\bf C}) \rightarrow {K_i}( {C^{*}}(\Gamma) )
\end{equation}
 which is known in many cases to be an isomorphism. 
 Here $B \Gamma$ is the classifying space of $\Gamma$. 
 Now for $\Gamma = {\bf Z}^3$, 
we have $B {\bf Z}^3 = {\bf T}^3$, hence the assembly map gives us an isomorphism 
\begin{equation}
\mu : {KK^i}( C( {\bf T}^3), {\bf C}) \cong {K_i}( C({\bf T}^3 ) )
\end{equation}
 which we can now describe explicitly. 
 We have 
\begin{equation*}
 \mu : {\bf z}_0 \mapsto \pm [1],
\end{equation*}
\begin{equation*}
{\psi_i^{*}}( {\bf w}_1) \mapsto \pm [U_i ],
\end{equation*}
\begin{equation*}
{\phi_i^{*}}({\bf Dirac} ) \mapsto \pm [ {Bott}_i],
\end{equation*}
\begin{equation}
{\bf Dirac} \mapsto \pm [ Bott].
\end{equation}
 This a very instructive example, in part because it is completely transparent. 
Note that although we found explicitly the eight different Fredholm modules generating the K-homology, only one of them, ${\bf Dirac}$, depends on the entire algebra, and this represents the fundamental class in K-homology. 
 
In the case of ${C^{*}}(H_3)$ the Fredholm module ${\partial_0} ( {\bf Dirac})$ is our candidate for the fundamental class of this noncommutative manifold. 
 If we were to try to build a three dimensional noncommutative geometry over 
${C^{*}}(H_3)$, we would start with the corresponding spectral triple and try to find the required extra structures. 
 It is not clear to me whether in general this ``fundamental class'' should exist at all, or whether it should be unique.

\section{Classification of derivations}

In this section we study the derivations of an appropriate smooth subalgebra $A^{\infty}$ of the group C*-algebra $A$ of the discrete Heisenberg group $H_3$, and we classify the derivations of the group ring ${\bf C}H_3$. 
 The algebra  $A$ is  
 generated by unitaries $U$, $V$ and $W$, with $W$ central,  satisfying the relation
\begin{equation}
VU=WUV.
\end{equation}
 We will be concerned with the ``smooth subalgebra'' 
\begin{equation}
\label{hbgsmsubalg}
{A^{\infty}} =
\{ \Sigma {a_{p,q,r}} {U^p}{V^q}{W^r} \, \},
\end{equation}
where 
$\{{a_{p,q,r}}\} \in S({\bf Z}^3) $, 
the set of functions on ${\bf Z}^3$ vanishing at infinity faster than any polynomial in $p$, $q$ and $r$.
Since ${H_3}$ is a finitely generated discrete amenable group, 
of polynomial growth of order four with respect to the word length function $L$ relative to the generators $U$, $V$ and $W$,  
it therefore has property (RD).
 Hence by Jolissaint's theorem [Jo89]  
the corresponding smooth subalgebra ${H_L^{\infty}} (\Gamma)$  is closed under holomorphic functional calculus. 
Furthermore, for the word length function $L$ defined above, 
${H_L^{\infty}}( \Gamma )$ coincides with $A^{\infty}$, 
hence $A^{\infty}$ is closed under holomorphic functional calculus.   
 So the inclusion map $i : A^{\infty} \hookrightarrow A$ induces isomorphisms $i_{*} : K_i (A^{\infty}) \cong K_i (A)$ on K-theory. 

Inspired in part by work of Bratteli, Elliott, Goodman and Jorgensen [BEGJ89] on classifying derivations of the noncommutative tori $A_{\theta}$, 
we classify all derivations from the group ring ${\bf C}H_3$ to  $A^{\infty}$ to itself, where by ``derivation'' we mean a linear map  
$\partial$
 satisfying the Leibnitz rule
 $\partial(ab) = \partial(a) b + a \partial(b)$.
 We will denote the vector space of all derivations
 $\partial : {\bf C}H_3 \rightarrow A^{\infty}$ by $Der({\bf C}H_3)$, and derivations
$\partial : A^{\infty} \rightarrow A^{\infty}$ 
 by $Der(A^{\infty})$. 
We begin with some preliminary lemmas. 
\begin{Lemma}
Let 
$\partial : {A^{\infty}} \rightarrow {A^{\infty}} $
be a derivation. 
Then $\partial(W)$ is central.
\end{Lemma} 
\begin{proof}
For any $x \in {A^{\infty}} $, 
we have 
\begin{equation}
\partial(Wx) = 
\partial(W) x + W \partial(x) =
\partial(xW) =
\partial(x) W + x \partial(W)
\end{equation}
Since $W$ is central, it follows that
$ \partial(W) x = x \partial(W)$ for all $x$.
\end{proof}

Let $Z(A^{\infty})$ be the centre of $A^{\infty}$. 
The vector spaces $Der({\bf C}H_3)$ and  $Der(A^{\infty})$ are naturally  left $Z(A^{\infty})$-modules, 
i.e. if $\partial$ is a derivation, so is the map
$ a \mapsto z \partial(a)$ for any $z \in Z(A^{\infty})$. 

\begin{Lemma}
The centre $Z(A^{\infty})$ of ${A^{\infty}}$ is  the smooth subalgebra ${C^{*}}(W) \cap A^{\infty}$ generated by $W$.
\end{Lemma} 
\begin{proof} 
 Given a central element $x \in {A^{\infty}}$, 
$x =
\Sigma {\alpha_{p,q,r}} {U^p}{V^q}{W^r}$ 
for some 
$\{  {\alpha_{p,q,r}} \} \in S({\bf Z}^3)$, and  
\begin{equation}
xU =
\Sigma {\alpha_{p-1,q,r-q}} {U^p}{V^q}{W^r}, \,
Ux =
\Sigma {\alpha_{p-1,q,r}} {U^p}{V^q}{W^r}
\end{equation}
\begin{equation}
xV =
\Sigma {\alpha_{p,q-1,r}} {U^p}{V^q}{W^r}, \,
Vx =
\Sigma {\alpha_{p,q-1,r-p}} {U^p}{V^q}{W^r}
\end{equation}
Equating coefficients we find that
 ${\alpha_{p,q,r}} = {\alpha_{p,q,r-q}} = {\alpha_{p,q,r-p}}$, 
for all $p$, $q$, $r$. 
 We require 
$\{ {\alpha_{p,q,r}} \} \in S({\bf Z}^3)$, 
 so we must have 
${\alpha_{p,q,r}} = 0$, unless 
$p = q = 0$.
Hence 
$ x = \Sigma {\alpha_{0,0,r}} {W^r} $,
which proves the assertion.
 \end{proof}

\begin{Lemma}
\label{partialW} 
 For any derivation $\partial : {A^{\infty}} \rightarrow {A^{\infty}} $, we have $\partial(W)=0$.
\end{Lemma} 
\begin{proof} 
Suppose that $\partial \in Der(A^{\infty})$ satisfies 
\begin{equation}
\partial(U) = \Sigma {a_{p,q,r}} {U^p}{V^q}{W^r}, \quad 
\partial(V) = \Sigma {b_{p,q,r}} {U^p}{V^q}{W^r},
\end{equation}
where 
$\{ a_{p,q,r} \}$, $\{ b_{p,q,r} \}$ are in $S({\bf Z}^3)$. 
Since $U$ and $V$ are unitary, it is immediate that 
\begin{equation}
\partial({U^{*}}) = - {U^{*}} \partial(U) {U^{*}}, \quad
\partial({V^{*}}) = - {V^{*}} \partial(V) {V^{*}}.
\end{equation}
 From the relation $W = VU{V^{*}}{U^{*}}$, 
we have 
\begin{equation}
\partial(W) = 
\partial(V) U{V^{*}}{U^{*}} +
V \partial(U) {V^{*}}{U^{*}} + 
VU \partial({V^{*}}) {U^{*}} +
VU{V^{*}} \partial({U^{*}}).
\end{equation}
Explicitly: 
\begin{equation}
\partial(W) =
\Sigma [ ( {b_{p,q+1,r-1}} - {b_{p,q+1,r+q-1}} )  +
({a_{p+1,q,r-p+q-1}} -{a_{p+1,q,r+q-1}}) ]
{U^p}{V^q}{W^r}. 
\end{equation}
The terms in $W^r$ (i.e. with $p=q=0$) are given by
\begin{equation}
{b_{0,1,r-1}} - {b_{0,1,r-1}} + {a_{1,0,r-1}} - {a_{1,0,r-1}} =0.
\end{equation}
So there are no terms in $W^r$ for any $r$, 
and since $\partial(W)$ must be central, it follows that $\partial(W)=0$. 
\end{proof}

Looking more closely at the relations above, we see that if we put $q=0$, we obtain
\begin{equation}
{a_{p+1,0,r-p-1}} - {a_{p+1,0,r-1}} = 0
\end{equation} 
for all $p$, $r$. 
It follows that  
${a_{p,0,r}}=0$ if $p \neq 1$.
In the same way, we see from putting $p=0$ that 
\begin{equation}
{b_{0,q+1,r-1}} - {b_{0,q+1,r+q-1}} = 0
\end{equation}
Hence ${b_{0,q,r}}=0$, unless $q=1$.
 In general we have the relation
\begin{equation}
( {b_{p,q+1,r-1}} - {b_{p,q+1,r+q-1}} )  +
({a_{p+1,q,r-p+q-1}} -{a_{p+1,q,r+q-1}}) = 0
\end{equation}
for all $p$, $q$ and $r$. 
This will be important in the sequel. 
 \begin{Defn}
We define derivations ${\partial_1}$, ${\partial_2}$ on ${\bf C}H_3$ 
by 
\begin{equation}
\label{hbgpartialone}
{\partial_1} (U) = U, \quad
{\partial_1} (V) = 0 , \quad
{\partial_1} (W) = 0,
\end{equation}
\begin{equation}
\label{hbgpartialtwo}
{\partial_2} (U) = 0, \quad
{\partial_2} (V) = V , \quad
{\partial_2} (W) = 0.
\end{equation}
\end{Defn} 
 We can extend $\partial_1$ and $\partial_2$ to $A^{\infty}$ via 
\begin{equation}
\partial_i ( \Sigma a_{p,q,r} U^p V^q W^r ) =
\Sigma a_{p,q,r} \partial_i (U^p V^q W^r) \quad (i=1,2).
\end{equation}
It follows from our previous work that these derivations are well-defined. 
 Note that since for each $n$, $\partial_1 (U^n) = n {U^{n-1}} \partial_1 (U)$,  $\partial_1$ (and similarly $\partial_2$) is an unbounded derivation and hence not inner. 

\begin{Defn} 
Given $x \in {A^{\infty}}$, 
$x =
\Sigma {\alpha_{p,q,r}} {U^p}{V^q}{W^r}$, 
 we denote by  
${\partial_{x}}$ 
the inner derivation corresponding to $x$ :  
\begin{equation}
 {\partial_{x}} (a) = [a,x] = ax-xa
\end{equation}
\end{Defn} 
 Note that, taking 
$x = \Sigma {\alpha_{p,q,r}} {U^p}{V^q}{W^r}$,
 we can choose ${\alpha_{0,0,r}}=0$.
From our previous calculations we see that 
\begin{equation}
{\partial_{x}}(U) = Ux-xU =
\Sigma ( {\alpha_{p-1,q,r}} - {\alpha_{p-1,q,r-q}} ) {U^p}{V^q}{W^r},
\end{equation}
\begin{equation}
{\partial_{x}}(V) = Vx-xV =
\Sigma ( {\alpha_{p,q-1,r-p}} - {\alpha_{p,q-1,r}} ) {U^p}{V^q}{W^r}.
\end{equation}
Note that in ${\partial_{x}}(U)$, all terms in ${U^p}{W^r}$ (i.e. with $q=0$) vanish, 
while similarly for ${\partial_{x}}(V)$, 
all terms in ${V^q}{W^r}$ (i.e. with $p=0$) are zero.
 Now we can state the main result: 
\begin{Thm}
\label{derivhbg}
As a left $Z(A^{\infty})$-module, 
$Der({\bf C}H_3)$ is spanned by the derivations $\partial_1$ and $\partial_2$, together with the inner derivations $\partial_x$.
Any derivation $\partial : {\bf C}H_3  \rightarrow {A^{\infty}} $
can be written uniquely as 
$\partial = {z_1} {\partial_1} + {z_2}{\partial_2} + {\partial_x}$, 
 for some $z_1$, $z_2 \in Z(A^{\infty})$, and $x \in A^{\infty}$.
\end{Thm} 
\begin{proof} 
Given 
$\partial : {\bf C}H_3  \rightarrow {A^{\infty}} $, 
defined by 
\begin{equation}
\partial(U) = \Sigma {a_{p,q,r}} {U^p}{V^q}{W^r}, \quad
\partial(V) = \Sigma {b_{p,q,r}} {U^p}{V^q}{W^r}, \quad
\partial(W) = 0,
\end{equation}
we will prove that 
\begin{equation}
\partial = 
({\Sigma_{r \in {\bf Z}}} {a_{1,0,r}} {W^r}){\partial_1} +
({\Sigma_{r \in {\bf Z}}} {b_{0,1,r} {W^r}}){\partial_2} +
{\partial_{x} }
\end{equation}
for an appropriate $x \in {A^{\infty}}$.
The infinite sums of derivations are interpreted in the following sense. 
We give the vector space of all derivations the weak topology in which a sequence of derivations ${ \{ {\partial_k} \} }_{k \in {\bf Z}}$ converges to a derivation $\partial$ if and only if 
${\partial_k}(a) \rightarrow \partial(a)$ 
(in norm) for all $a \in {A^{\infty}}$. 
Choosing 
\begin{equation}
{x =
\Sigma \, {\alpha_{p,q,r}} {U^p}{V^q}{W^r}},
\end{equation}
we have
\begin{equation}
\partial_{x}(U) = \Sigma ( {\alpha_{p-1,q,r}} - {\alpha_{p-1,q,r-q}} ) {U^p}{V^q}{W^r},
\end{equation}
\begin{equation}
\partial_{x}(V) = \Sigma ( {\alpha_{p,q-1,r-p}} - {\alpha_{p,q-1,r}} ) {U^p}{V^q}{W^r}.
\end{equation}
Hence we need to find 
$\{ {\alpha_{p,q,r}} \} \in S({\bf Z}^3) $ 
such that, for $p \neq 0$, $q \neq 0$ 
\begin{equation}
{\alpha_{p-1,q,r}} - {\alpha_{p-1,q,r-q}} = a_{p,q,r},
\end{equation}
\begin{equation}
{\alpha_{p,q-1,r-p}} - {\alpha_{p,q-1,r}} = b_{p,q,r}.
\end{equation}

 We do this explicitly. 
 We need to solve these equations ensuring that our solutions 
$\{ {\alpha_{p,q,r}} \}$ 
 are Schwarz functions on ${\bf Z}^3$. 
There are eight different cases to be considered, corresponding to the different possibilities for the signs of $p$, $q$, $r$.

{\bf Case 1:} (+\,+\,+) $p, q > 0$, $r \geq 0$. Define 
\begin{equation}
\label{sumone}
{\alpha_{p,q,r}} = - {\Sigma_{k=1}^{\infty}} {a_{p+1,q,r+kq}},
\end{equation} 
\begin{equation}
\label{sumtwo}
{\alpha_{p,q,r}} = {\Sigma_{k=1}^{\infty}} {b_{p,q+1,r+kp}}.
\end{equation} 
We verify all the details for this case. 
Since $\{ {a_{p,q,r}} \}$ and 
$\{ {b_{p,q,r}}\}$ are Schwarz functions, 
it is immediate that the infinite sums converge absolutely. 
We first of all need to check that the two infinite sums (110), (111) are in fact equal to one another. 
Since
\begin{equation*}
{a_{p+1,q,r+q}} = {a_{p+1,q,r+q+p}} - {b_{p,q+1,r+p}} + {b_{p,q+1,r+q+p}},
\end{equation*}
\begin{equation*}
{a_{p+1,q,r+q+p}} = {a_{p+1,q,r+q+2p}} - {b_{p,q+1,r+2p}} + {b_{p,q+1,r+q+2p}},
\end{equation*} 
and so on, it follows that
\begin{equation}
{a_{p+1,q,r+q}} = -{\Sigma_{k=1}^{\infty}}{b_{p,q+1,r+kp}}  +  {\Sigma_{k=1}^{\infty}}{b_{p,q+1,r+q+kp}},
\end{equation}
and therefore
\begin{equation}
{a_{p+1,q,r+jq}} = -{\Sigma_{k=1}^{\infty}}{b_{p,q+1,r+(j-1)q+kp}}  +  {\Sigma_{k=1}^{\infty}}{b_{p,q+1,r+jq+kp}}.
\end{equation}
Therefore
\begin{equation}
{\Sigma_{j=1}^{\infty}}{a_{p+1,q,r+jq}} = -{\Sigma_{j=0,k=1}^{\infty}}{b_{p,q+1,r+jq+kp}}  +  {\Sigma_{j=1,k=1}^{\infty}}{b_{p,q+1,r+jq+kp}} 
= - {\Sigma_{k=1}^{\infty}}{b_{p,q+1,r+kp}},
\end{equation}
which proves the equality of our two expressions for 
${\alpha_{p,q,r}}$. 
We also need to verify the relations  
\begin{equation}
{\alpha_{p-1,q,r}} - {\alpha_{p-1,q,r-q}} = a_{p,q,r}, \quad
{\alpha_{p,q-1,r-p}} - {\alpha_{p,q-1,r}} = b_{p,q,r}.
\end{equation}
This is easily done: 
\begin{equation}
{\alpha_{p-1,q,r}} = -{\Sigma_{k=1}^{\infty}} {a_{p,q,r+kq}},
\end{equation}
while
\begin{equation}
{\alpha_{p-1,q,r-q}} = -{\Sigma_{k=0}^{\infty}} {a_{p,q,r+kq}}.
\end{equation}
Furthermore 
\begin{equation}
{\alpha_{p,q-1,r-p}} = {\Sigma_{k=0}^{\infty}} {b_{p,q,r+kp}},
\end{equation}
\begin{equation}
{\alpha_{p,q-1,r}} = {\Sigma_{k=1}^{\infty}} {b_{p,q,r+kp}}.
\end{equation}

The remaining seven cases are as follows. In each case we give an explicit formula for the $a_{p,q,r}$ as two different infinite sums, which exactly as above, can be shown to be equal, and to satisfy all the desired relations. 

{\bf Case 2:} (+\,+\,-) $p,q > 0$, $r \leq 0$.  Define 
\begin{equation}
{\alpha_{p,q,r}} = {\Sigma_{k=0}^{\infty}} {a_{p+1,q,r-kq}},
\end{equation} 
\begin{equation}
{\alpha_{p,q,r}} = -{\Sigma_{k=0}^{\infty}} {b_{p,q+1,r-kp}}.
\end{equation} 

{\bf Case 3:} (+\,-\,+) $p>0$, $q<0$, $r \geq 0$.  Define 
\begin{equation}
{\alpha_{p,q,r}} = {\Sigma_{k=0}^{\infty}} {a_{p+1,q,r-kq}},
\end{equation} 
\begin{equation}
{\alpha_{p,q,r}} = {\Sigma_{k=1}^{\infty}} {b_{p,q+1,r+kp}}.
\end{equation} 

{\bf Case 4:} (+\,-\,-) $p>0$, $q<0$, $r \leq 0$. Define 
\begin{equation}
{\alpha_{p,q,r}} = - {\Sigma_{k=1}^{\infty}} {a_{p+1,q,r+kq}},
\end{equation} 
\begin{equation}
{\alpha_{p,q,r}} = - {\Sigma_{k=0}^{\infty}} {b_{p,q+1,r-kp}}.
\end{equation} 

{\bf Case 5:} (-\,+\,+) $p <0$, $q>0$, $r \geq 0$. Define 
\begin{equation}
{\alpha_{p,q,r}} = - {\Sigma_{k=1}^{\infty}} {a_{p+1,q,r+kq}},
\end{equation} 
\begin{equation}
{\alpha_{p,q,r}} = -{\Sigma_{k=0}^{\infty}} {b_{p,q+1,r-kp}}.
\end{equation} 

{\bf Case 6:} (-\,+\,-) $p<0$, $q>0$, $r \leq 0$.  Define 
\begin{equation}
{\alpha_{p,q,r}} =  {\Sigma_{k=0}^{\infty}} {a_{p+1,q,r-kq}},
\end{equation} 
\begin{equation}
{\alpha_{p,q,r}} = {\Sigma_{k=1}^{\infty}} {b_{p,q+1,r+kp}}.
\end{equation} 

{\bf Case 7:} (-\,-\,+) $p, q <0$, $r \geq 0$.  Define 
\begin{equation}
{\alpha_{p,q,r}} =  {\Sigma_{k=0}^{\infty}} {a_{p+1,q,r-kq}},
\end{equation} 
\begin{equation}
{\alpha_{p,q,r}} = - {\Sigma_{k=0}^{\infty}} {b_{p,q+1,r-kp}}.
\end{equation} 

{\bf Case 8:} (-\,-\,-) $p, q < 0$, $r \leq 0$. Define 
\begin{equation}
{\alpha_{p,q,r}} = - {\Sigma_{k=1}^{\infty}} {a_{p+1,q,r+kq}},
\end{equation} 
\begin{equation}
{\alpha_{p,q,r}} = {\Sigma_{k=1}^{\infty}} {b_{p,q+1,r+kp}}.
\end{equation} 

This completes the proof. 
\end{proof}

This gives a very complete and satisfactory picture of the derivations of $A^{\infty}$.
What is surprising is the result (Lemma 27) that $\partial(W) = 0$ for every $\partial \in Der(A^{\infty})$. 
In the picture of $A = {C^{*}}(H_3)$ as a continuous field of C*-algebras over the circle, this tells us that derivations only act fiberwise - there is no derivation corresponding to ``going round the circle''. 
 From the viewpoint of derivations, this group C*-algebra does not really look three-dimensional.

\section{Cyclic cohomology of the group ring}

We use  Burghelea's theorem [Bu85]  to calculate the cyclic cohomology of the group ring of the discrete Heisenberg group $H_3$. 
 Recall that for a countable discrete group $G$, 
Burghelea's theorem states that the cyclic cohomology of the group ring ${\bf C}G$ is given [Co94], p213,  by 
\begin{equation}
{HC^{*}}({\bf C}G) = {\Pi_{ {\hat{g}} \in <G>'} ( {H^{*}}( {N_g};{\bf C}) \otimes {HC^{*}}({\bf C}) )} \times  {\Pi_{ {\hat{g}} \in <G>''} {H^{*}}( {N_g}; {\bf C})}.
\end{equation}
  Given $g \in G$, we let ${C_g} = \{ h \in G : gh=hg \} $ be the centraliser of $g$, and 
 ${N_g} = {C_g}/ {g^{\bf Z}}$ be the quotient of ${C_g}$ by the central normal subgroup generated by $g$. 
Here $<G>$ denotes the set of conjugacy classes of elements of $G$, while $<G>'$ is the set of conjugacy classes of elements of finite order, and $<G>''$ is the set of conjugacy classes of elements of infinite order. 

For $G= H_3$, the only element of finite order is the identity $e$, and $C_e = N_e = H_3$. 
For any other element $x = {U^p}{V^q}{W^r}$, we have ${U^p}{V^q}{W^r} \sim {U^p}{V^q}{W^{r+nk}}$, 
where $k$ is the highest common factor of $p$ and $q$, and $n$ is an integer. 

\begin{Prop} 
\label{cxnxhbg}
For each $g \in {H_3}$, $g \neq e$, we have 
$C_g \cong {\bf Z}^2$, unless $g \in Z(H_3)$, in which case $C_g \cong H_3$. 
In the first case  $N_g ={\bf Z} \times {\bf Z}_l$, while in the second 
$N_g$ is nonabelian, being a central extension of  ${\bf Z}^2$ by 
 ${\bf Z}_l$ 
for some $l \in {\bf Z}_{+}$ (depending on $g$) that we will determine.
\end{Prop}
\begin{proof}
Suppose that 
$g = 
\left(
\begin{array}{ccc}
1 & p & r \cr
0 & 1 & q \cr
0 & 0 & 1 \cr
\end{array}
\right)
$. 
 Then there are four separate cases to be considered. 
 Note that each case is stable under conjugation. 

{\bf Case 1:} $p$, $q$ are both nonzero. 
 Let $k = hcf(p,q)$ be the highest common factor of $p$ and $q$, 
 and let $p = kp'$, $q = kq'$, where now $hcf(p',q') =1$.  
Then an easy calculation shows that 
\begin{equation}
C_g = 
 \{
\left(
\begin{array}{ccc}
1 & np' & c \cr
0 & 1 & nq' \cr
0 & 0 & 1 \cr
\end{array}
\right)
 \quad : \quad
 n, c \in {\bf Z}
\}
\end{equation}
 It is easy to check that $C_g$ is abelian, and in fact 
 $C_g \cong {\bf Z}^2$ via the isomorphism 
\begin{equation}
\left(
\begin{array}{ccc}
1 & np' & c \cr
0 & 1 & nq' \cr
0 & 0 & 1 \cr
\end{array}
\right)
 \mapsto
 (n, c - S_n)
\end{equation}
 where for each $n \in {\bf Z}$, we define 
 $S_n = {\frac{1}{2}}p'q'n(n-1)$. 
 Under this isomorphism, 
we have $g \mapsto (k, r - {S_k})$, 
and we can identify $N_g$ with the quotient of ${\bf Z}^2$ 
by the normal subgroup (isomorphic to ${\bf Z}$) 
generated by this element.
 The following Lemma we leave as an elementary exercise for the reader. 

\begin{Lemma}
\label{eltlemma}
 Given $(a,b) \in {\bf Z}^2$, $(a,b) \neq (0,0)$, let $G$ be the subgroup (isomorphic to ${\bf Z}$)  generated by this element. Then ${\bf Z}^2 / G \cong {\bf Z} \times {\bf Z}_l$, where $l = hcf(a,b)$. If either $a$ or $b$ is zero, we take $l = max( |a|, |b|)$.
\end{Lemma}

Hence ${N_g} \cong {\bf Z} \times {\bf Z}_l$, where 
 $l = hcf(k, r -S_k)$.

{\bf Case 2:} $p \neq 0$, $q=0$.  In this case 
\begin{equation}
C_g \cong \{
\left(
\begin{array}{ccc}
1 & a & c \cr
0 & 1 & 0 \cr
0 & 0 & 1 \cr
\end{array}
\right)
\}
 \cong {\bf Z}^2.
\end{equation}
 Hence by Lemma 32 we have $N_g \cong {\bf Z} \times {\bf Z}_l$, where $l = hcf(p,r)$.

{\bf Case 3:} $p =0$, $q \neq 0$. This proceeds exactly as in Case 2, to give 
 $C_g \cong {\bf Z}^2$, and $N_g \cong {\bf Z} \times {\bf Z}_l$, where $l = hcf(q,r)$.

{\bf Case 4a:} $p=0$, $q=0$, $r= \pm 1$.
 Then $g$ generates $Z(H_3)$, and $N_g \cong$  ${\bf Z}^2$. 

{\bf Case 4b:} $p=0$, $q=0$ (and $r \neq 0$, $\pm 1$).
 So $g \in Z(H_3)$, hence $C_g = H_3$. 
 In this case $N_g$ is not abelian, but it is a central extension
\begin{equation}
0 \rightarrow {\bf Z}_{|r|} 
\rightarrow N_g 
\rightarrow {\bf Z}^2
\rightarrow 0
\end{equation}
 which we will use to calculate the group cohomology.
 This completes the proof of Prop 31.
\end{proof}

Now, we know from [Wei94], p166, that (for $l \neq 1$) 
\begin{equation}
{H^n}( {\bf Z} \times {\bf Z}_l ; {\bf C} ) \cong 
{\oplus_{p=0}^{n}} {H^p}( {\bf Z} ; {\bf C}) \otimes {H^{n-p}}( {\bf Z}_l ; {\bf C}) 
\end{equation}
with 
\begin{equation}
 {H^p}( {\bf Z}; {\bf C}) \cong
\left\{
\begin{array}{cc}
{\bf C} & : p=0,1 \cr
0 & : p \geq 2 \cr
\end{array}
\right.
\end{equation}
\begin{equation}
 {H^q}( {\bf Z}_l; {\bf C}) \cong
\left\{
\begin{array}{cc}
{\bf C} & : q=0 \cr
0 & : q \geq 1 \cr
\end{array}
\right.
\end{equation}
 (for $l \geq 2$),  
so we have 
\begin{equation}
{H^n}( {\bf Z} \times {\bf Z}_l ; {\bf C} ) \cong 
\left\{ 
\begin{array}{cc}
{\bf C} & : {n=0, 1} \\
 0 & : {n \geq 2} \\
\end{array}
\right.
\end{equation}
 Obviously, if $l=1$, we have 
${\bf Z} \times {\bf Z}_l \cong {\bf Z}$.
 So this deals with Cases 1, 2 and 3.  Now,
\begin{equation}
{H^p}( {\bf Z}^2 ; {\bf C}) \cong 
 {\oplus}_{q=0}^{p}
 {H^q}({\bf Z} ; {\bf C}) \otimes 
{H^{p-q}}( {\bf Z} ; {\bf C})
\cong 
\left\{
\begin{array}{cc}
{\bf C} & : p=0 \cr
{\bf C}^2 & : p=1 \cr
{\bf C} & : p=2 \cr
0 & : p \geq 3 \cr
\end{array}
\right.
\end{equation}
 This takes care of Case 4a. 

For Case 4b, we use the Hochschild-Serre spectral sequence [Wei94], p195. 
 For the central extension 
 \begin{equation}
0 \rightarrow {\bf Z}_{|r|} 
\rightarrow N_g 
\rightarrow {\bf Z}^2
\rightarrow 0
\end{equation}
 we have [Wei94], p196, that 
${\bf Z}^2$ acts trivially on 
$H^n ({\bf Z}_{|r|} ; {\bf C})$. 
 Define 
\begin{equation}
E_2^{p,q} =
H^p ( {\bf Z}^2 ; H^q ({\bf Z}_{|r|} ; {\bf C})). 
\end{equation}
 Then $E_2^{p,q}$ converges to
$H^{p+q}( N_g ; {\bf C})$. 
We have 
\begin{equation}
E_2^{p,q} =
\left\{
\begin{array}{cc}
H^p ( {\bf Z}^2 ; {\bf C}) & : q=0 \cr
0 & : q \neq 0 \cr
\end{array}
\right.
\end{equation}
So
\begin{equation}
E_2^{p,0} =
\left\{
\begin{array}{cc}
{\bf C} & : p=0 \cr
{\bf C}^2 & : p =1 \cr
{\bf C} & : p=2 \cr
0 & : p \geq 3
\end{array}
\right.
\end{equation}
 which implies that
\begin{equation}
H^n (N_g ; {\bf C}) =
\left\{
\begin{array}{cc}
{\bf C} & : n=0 \cr
{\bf C}^2 & : n =1 \cr
{\bf C} & : n=2 \cr
0 & : n \geq 3
\end{array}
\right.
\end{equation}

If $g \in Z(H_3)$, then $g$ is conjugate only to itself, hence there are countably many conjugacy classes of 
 $\hat{g} \in <G>''$, for which $g \in Z(H_3)$.
 Similarly there are countably many 
conjugacy classes of 
 $\hat{g} \in <G>''$, for which $g$ is not an element of $Z(H_3)$.
Therefore, 
\begin{equation}
{\Pi_{ {\hat{g}} \in <G>''} {H^{*}}( {N_g}; {\bf C})}
 \cong 
\left\{ 
\begin{array}{cc}
{\Pi_{\bf Z}}{\bf C} & : {n=0,1,2} \\
 0 & : {n \geq 3} \\
\end{array}
\right.
\end{equation}

We also need the following :

\begin{Prop} The group cohomology of ${H_3}$ is given by  
\begin{equation}
{H^n}( {H_3} ; {\bf C} ) = 
\left\{ 
\begin{array}{cc}
{\bf C} & : {n=0} \\
{\bf C}^2  &   : {n =1} \\
{\bf C}^2 & : {n=2} \\
{\bf C} & : {n=3} \\
0 & : {n > 3}\\
\end{array}
\right.
\end{equation}
\end{Prop} 
 \begin{proof}
 We note first of all that since ${\bf C}$ is a trivial ${H_3}$-module, ${H^0}( {H_3} ; {\bf C} ) = {\bf C}$. 
Furthermore, we know that 
\begin{equation}
{H^1}( {H_3} ; {\bf C} ) \cong {Hom_{Groups}}( {H_3} , {\bf C} ) \cong {Hom_{Ab}}( {H_3}/ [ {H_3},{H_3}] , {\bf C}) \cong {\bf C}^2 
\end{equation} 
 
Since  $H_3$ is a central extension 
\begin{equation}
0 \rightarrow {\bf Z} \rightarrow {H_3} \rightarrow {\bf Z}^2 \rightarrow 0 
\end{equation}
 we can use the Hochschild-Serre spectral sequence to calculate the rest of the group cohomology.  
 Since here ${\bf Z}$ is the centre of ${H_3}$, the ${\bf Z}^2$ action on ${H^n}( {\bf Z} ; {\bf C}) $ is trivial. 
 We have 
${E_{2}^{p,q}} = {H^p}( {\bf Z}^2 ; {H^q}( {\bf Z} ; {\bf C}) )$. 
 Hence ${E_{2}^{p,q}} =0$ unless $p=0,1,2$ and $q=0,1$. 
We have 
\begin{equation}
{E_{2}^{0,0}} = {\bf C}, \, {E_{2}^{0,1}} = {\bf C} 
\end{equation}
\begin{equation}
{E_{2}^{1,0}} = {\bf C}^2, \, {E_{2}^{1,1}} = {\bf C}^2 
\end{equation}
\begin{equation}
{E_{2}^{2,0}} = {\bf C}, \, {E_{2}^{2,1}} = {\bf C} 
\end{equation}
 The only potentially nontrivial differential is 
$d : {E_{2}^{0,1}} \rightarrow {E_{2}^{2,0}}$ 
but we already know that 
 ${H^1}( {H_3} ; {\bf C} )$  $\cong {\bf C}^2 $, 
 hence $d =0$. 
So everything is over at the first step of the spectral sequence.  
\end{proof}

Combining these results, we have :

\begin{Thm} 
\label{cchbggpring}
The cyclic cohomology of the group ring ${\bf C}{H_3}$ of the discrete Heisenberg group ${H_3}$ is given by 
\begin{equation}
{HC^n}({\bf C}{H_3}) = 
\left\{ 
\begin{array}{cc}
{\Pi}_{\bf Z} {\bf C} & : {n=0,1,2} \\
{\bf C}^3 & : {n \geq 3} \\
\end{array}
\right.
\end{equation}
\end{Thm}

\begin{Cor}
 Both the even and the odd periodic cyclic cohomology of ${\bf C}{H_3}$ are isomorphic to ${\bf C}^3$. 
\end{Cor}

Exhibiting the generating cyclic cocycles, and calculating their pairings with the generators of K-theory, is involved, and we will not present these calculations here. 

We note finally that the cyclic cohomology of the ``smooth subalgebra'' 
$A^{\infty}$ defined in (85) can be calculated using results of Nest [Ne88] on smooth crossed products by ${\bf Z}$. 

\begin{Prop}
The cyclic cohomology of the smooth subalgebra $A^{\infty}$ is given by 
${HC^n}(A^{\infty}) = {\bf C}^3$,  $n \geq 2$. 
 In particular, the even and odd periodic cyclic cohomology are both isomorphic to ${\bf C}^3$. 
\end{Prop}

\section{Acknowledgements}
 I would like to thank my advisor Professor Marc Rieffel for his advice and support throughout my time in Berkeley. I am extremely grateful for his help. 
 I would also like to thank Erik Guentner, Nate Brown and Frederic Latremoliere for many very useful discussions.

\bibliographystyle{plain}

\end{document}